\documentclass{icmart}

\usepackage{mathrsfs}
\usepackage[dvipsnames]{xcolor}
\usepackage[colorlinks=true,allcolors=BrickRed]{hyperref}
\usepackage{enumitem}
\usepackage{breakurl}
\contact[rdlyons@indiana.edu]{Dept. of Math., Indiana University, 831 E.
3rd St., Bloomington, IN 47405-7106 USA}

\numberwithin{equation}{section}

\theoremstyle{definition}

\def\thmenv#1#2#3{\begin{#1} \label{#1:#2} #3 \end{#1}}
\def\richthmenv#1#2#3#4{\begin{#1}[#3] \label{#1:#2} #4 \end{#1}}

\def\procl#1.#2 #3\endprocl{%
       \ifx#1t\thmenv{thm}{#2}{#3}\fi
       \ifx#1l\thmenv{lem}{#2}{#3}\fi
       \ifx#1p\thmenv{pro}{#2}{#3}\fi
       \ifx#1c\thmenv{cor}{#2}{#3}\fi
       \ifx#1d\thmenv{dfn}{#2}{#3}\fi
       \ifx#1g\thmenv{conj}{#2}{#3}\fi
       \ifx#1q\thmenv{question}{#2}{#3}\fi
       \ifx#1r\thmenv{remark}{#2}{{\rm #3}}\fi
    }%

\def\rprocl#1.#2 #3 #4\endprocl{%
       \ifx#1t\richthmenv{thm}{#2}{#3}{#4}\fi
       \ifx#1l\richthmenv{lem}{#2}{#3}{#4}\fi
       \ifx#1p\richthmenv{pro}{#2}{#3}{#4}\fi
       \ifx#1c\richthmenv{cor}{#2}{#3}{#4}\fi
       \ifx#1d\richthmenv{dfn}{#2}{#3}{#4}\fi
       \ifx#1g\richthmenv{conj}{#2}{#3}{#4}\fi
       \ifx#1q\richthmenv{question}{#2}{#3}{#4}\fi
       \ifx#1r\richthmenv{remark}{#2}{#3}{{\rm #4}}\fi
    }%

\def\rref#1.#2/{%
      \ifx #1sSection~\ref{s.#2}\fi
      \ifx #1SSubsection~\ref{S.#2}\fi
      \ifx #1tTheorem~\ref{thm:#2}\fi  
      \ifx #1lLemma~\ref{lem:#2}\fi 
      \ifx #1cCorollary~\ref{cor:#2}\fi 
      \ifx #1pProposition~\ref{pro:#2}\fi 
      \ifx #1dDefinition~\ref{dfn:#2}\fi
      \ifx #1gConjecture~\ref{conj:#2}\fi 
      \ifx #1qQuestion~\ref{question:#2}\fi 
      \ifx #1rRemark~\ref{remark:#2}\fi 
      \ifx #1aAppendix~\ref{a.#2}\fi 
      \ifx #1fFigure~\ref{f.#2}\fi
      \ifx #1e(\ref{e.#2})\fi
      \ifx #1b\cite{#2}\fi
        }

\def\rlabel #1 #2{\begin{equation} \label{#1} #2 \end{equation}}

\def\rproof{\begin{proof}}

\newenvironment{proofof}[1]{{\medbreak\noindent \em Proof of #1.  }}{\hfill\qed\medbreak}

\def\Qed{\end{proof}}

\def\bsection#1#2{\bigbreak\section{#1}\label{#2}}
\def\bsubsection#1#2{\bigbreak\subsection{#1}\label{#2}}

\def\bbC{{\mathbb C}}

\def\bbZ{{\mathbb Z}}
\def\bbN{{\mathbb N}}

\def\bbR{{\mathbb R}}

\def\bbT{{\mathbb T}}

\def\bfone{{\boldsymbol 1}}

\def\verts{\mathsf{V}}
\def\vertex{\mathsf{V}}
\def\edges{\mathsf{E}}
\def\edge{\mathsf{E}}

\def\bE{{\bf E}}
\def\bP{{\bf P}}
\def\Seq#1{\langle #1 \rangle}

\def\dfnterm#1{\textit{\textbf{#1}}}

\def\st{\,;\;}
\def\arXiv#1{\url{http://www.arxiv.org/abs/#1}}
\def\restrict{\mathord{\upharpoonright}}   %% for restriction
\def\gp{{\Gamma}}
\def\gpe{{\gamma}}
\def\dom{\preccurlyeq}   %% stochastically dominated by
\font\frak=eufm10   %% or  eufb10
\font\scriptfrak=eufm7
\font\scriptscriptfrak=eufm5
\def\mathfrak#1#2{%       %% This cannot be done as
\def#1{{\mathchoice%
{{\hbox{\frak #2}}}%
{{\hbox{\frak #2}}}%
{{\hbox{\scriptfrak #2}}}%
{{\hbox{\scriptscriptfrak #2}}}}}}
\mathfrak{\ba}{B}  %% random variable for base
\mathfrak{\qba}{S}  %% random variable for subset
\mathfrak{\fo}{F}  %% random forest
\mathfrak{\rH}{H}  %% random subspace
\mathfrak{\rQ}{Q}  %% random positive contraction
\mathfrak{\rv}{v}  %% random unit vector
\mathfrak{\rM}{M}  %% random matrix
\def\fsf{\mathsf{FSF}}         %% Weighted free spanning forest measure
\def\wsf{\mathsf{WSF}}         %% Weighted wired spanning forest measure
\def\ust{\mathsf{UST}}         %% uniform spanning tree measure
\def\STAR{\bigstar}        %% The space spanned by the stars.
\def\CYCLE{\diamondsuit}
\def\cbuldot{{\raise.25ex\hbox{$\scriptscriptstyle\bullet$}}}
\def\bp{o}
\def\boI#1{{\bf 1}_{#1}}
\def\II#1{\boI{\{#1\}}}
\def\dbar{\bar d}

\def\EBig#1{\bE\Big[#1\Big]}
\def\Ebig#1{\bE\big[#1\big]}
\def\Var{{\rm Var}}
\def\BLPSusf{\cite{BLPS:usf}}
\def\ip#1{(\changecomma #1)}
\def\bigip#1{\big(\bigchangecomma #1\big)}
\def\Bigip#1{\Big(\bigchangecomma #1\Big)}

\def\leftip#1{\left(\leftchangecomma #1\right)}  
\def\changecomma#1,{#1,\,}
\def\bigchangecomma#1,{#1,\;}
\def\leftchangecomma#1,{#1,\ }
\def\sgn(#1){(-1)^{#1}}
\def\eH{{\ell^2(E)}}
\def\Ext{{\rm Ext}}
\def\bu{{\bf u}}
\def\bv{{\bf v}}
\def\bw{{\bf w}}
\def\all#1{\forall #1\enspace}

\def\cN{{\mathcal N}}

\def\sF{{\mathscr F}}
\mathfrak{\sX}{X}
\def\cA{{\mathcal A}}
\def\ev#1{{\mathcal #1}}
\def\ent{\mathsf {Ent}}  %% entropy
\def\asp{\alpha}  % aspect
\def\bm{\mathsf {BM}}  % B-M density
\def\Seq#1{\langle #1 \rangle}
\def\st{\,;\;}
\def\Bern{{\rm Bern}}
\def\Bin{{\rm Bin}}
\def\all#1{\forall #1\enspace}
\def\texists#1{\exists #1\enspace}
\def\diag{\Delta}  % diagonal of product space
\DeclareMathOperator{\img}{im}

\def\isc{\mathscr{U}}
  
\def\medstrut{\vrule height12pt depth5pt width0pt}   
\def\smallstrut{\vrule height10pt depth5pt width0pt}   
\def\poly{\mathsf{Poly}}
\DeclareMathOperator{\myIm}{Im}
\DeclareMathOperator{\myRe}{Re}
\def\disk{\mathbb {D}}
\def\GM{\mathsf{GM}}  %% geometric mean
\def\str{\star}          %% A star.
\def\gLoew{\succeq}
\def\lLoew{\preceq}
\def\FK{\mathsf{FK}}
\def\mv{\omega}  % multivector to represent the highest ext. power
\def\Sym{\mathsf{Sym}}  % symmetric group and extensions
\def\subsp{\le}  % closed subspace
\def\mud{\mu_{\rm d}}
\def\muc{\mu_{\rm c}}
\DeclareMathOperator{\tr}{tr}  % trace

\title[Determinantal Probability]{Determinantal Probability \\ 
\vspace{4mm}%
{\large Basic Properties and Conjectures}
}

\author[Russell Lyons]
{Russell Lyons\thanks{Partially supported by NSF grant DMS-1007244.}}

\begin{document}

\begin{abstract}
We describe the fundamental constructions and properties of determinantal
probability measures and point processes, giving streamlined proofs.
We illustrate these with some important examples.
We pose several general questions and conjectures.
\end{abstract}

\begin{classification}
Primary 
60K99, % None of the above, but in this section
60G55; % Point processes
Secondary 
42C30, % Completeness of sets of functions
37A15, % General groups of measure-preserving transformations
37A35, % Entropy and other invariants, isomorphism, classification
37A50, % Relations with probability theory and stochastic processes
68U99. % None of the above, but in this section
\end{classification}

\begin{keywords}
Random matrices, eigenvalues, orthogonal projections, positive
contractions, exterior algebra, stochastic domination, negative
association, point processes, mixtures, spanning trees, orthogonal
polynomials, completeness, Bernoulli processes. 
\end{keywords}

\maketitle

%\tableofcontents

\bsection{Introduction}{s.intro}

Determinantal point processes were
originally defined by Macchi \cite{Macchi} in physics.
Starting in the 1990s, determinantal probability
began to flourish as examples appeared in
numerous parts of mathematics \cite{Soshnikov:survey,Johansson,
Borodin:handbook}.
Recently, applications to machine learning have appeared \cite{KulTas}.

A discrete determinantal probability measure is one whose elementary
cylinder probabilities are given by determinants.
More specifically, suppose that $E$ is a finite or countable set and that
$Q$ is an $E \times E$ matrix.
For a subset $A \subseteq E$, let $Q\restrict A$ denote the submatrix of
$Q$ whose rows and columns are indexed by $A$.
If $\qba$ is a random subset of $E$ with the property that for all finite $A
\subseteq E$, we have 
\rlabel e.DPM
{\bP[A \subseteq \qba] = \det (Q\restrict A)
\,, }
then we call $\bP$ a \dfnterm{determinantal probability measure}.
The inclusion-exclusion principle in combination with \eqref{e.DPM}
yields the probability of each elementary
cylinder event.
Therefore,
for every $Q$, there is at most one probability
measure, to be denoted $\bP^Q$, on subsets of $E$ that satisfies \eqref{e.DPM}.
Conversely, it is known (see, e.g., \rref b.Lyons:det/) that there is a
determinantal probability measure corresponding to $Q$ if $Q$ is the matrix of
a positive contraction on $\ell^2(E)$ (in the standard orthonormal basis).

Technicalities are required even to define
the corresponding concept of determinantal point process for $E$
being Euclidean space or a more general space.
We present a virtually complete development of their basic properties
in a way that minimizes such technicalities by adapting the approach
of \rref b.Lyons:det/ from the discrete case.
In addition, we use an idea of Goldman \rref b.Goldman/
to deduce properties of the general case from corresponding properties in
the discrete case.

Space limitations prevent mention of most of what is known in
determinantal probability theory, which pertains largely to the analysis of
specific examples.
We focus instead
on some of the basic properties that hold for all determinantal processes
and on some intriguing open questions.

\bsection{Discrete Basics}{s.discrete}

Let $E$ be a denumerable set.

We identify a subset of $E$ with an element of $\{0, 1\}^E = 2^E$ in the
usual way.
There are several approaches to prove the basic existence results and
identities for determinantal probability measures. We sketch the one used
by \rref b.Lyons:det/.
This depends on understanding first the case where $Q$ is the matrix of an
orthogonal projection.
It also relies on exterior algebra so that the existence becomes immediate.

Any unit vector $v$ in a Hilbert space with orthonormal basis $E$ gives a
probability measure $\bP^v$ on $E$, namely, $\bP^v\big(\{e\}\big) :=
|\ip{v, e}|^2$ for $e \in E$. Applying this simple idea to multivectors
instead, we 
obtain the probability measures $\bP^H$ associated to orthogonal projections
$P_H$.
We refer to \rref b.Lyons:det/ for details not given here.

\bsubsection{Exterior Algebra}{S.ext}

Identify $E$ with the standard
orthonormal basis of the real or complex Hilbert space $\ell^2(E)$.
For $k \ge 1$, let $E_k$ denote a collection
of ordered $k$-element subsets of $E$ such that
each $k$-element subset of $E$ appears exactly once in $E_k$ in some
ordering. Define 
$$
\Lambda^k E := \bigwedge\nolimits^{\!k} \eH
:= \ell^2\Bigl( \big\{e_1 \wedge \cdots
\wedge e_k \st \Seq{e_1, \ldots, e_k} \in E_k \big\} \Bigr)
\,.
$$
If $k > |E|$, then $E_k = \emptyset$ and $\Lambda^k
E = \{0\}$.  We also define $\Lambda^0 E$ to be 
the scalar field, $\bbR$ or $\bbC$.
The elements of $\Lambda^k E$ are called \dfnterm{multivectors} of
\dfnterm{rank} $k$, or \dfnterm{$k$-vectors} for short. We then
define the \dfnterm{exterior} (or \dfnterm{wedge}) \dfnterm{product}
of multivectors in the usual alternating
multilinear way: $\bigwedge_{i=1}^k e_{\sigma(i)} = \sgn(\sigma)
\bigwedge_{i=1}^k e_i$ for any permutation $\sigma \in \Sym(k)$,
%of $\{1, 2, \ldots, k\}$ 
and $$\bigwedge_{i=1}^k
\sum_{e\in E'} a_i(e) e = \sum_{e_1, \ldots, e_k \in E'} \prod_{j=1}^k
a_j(e_j) \bigwedge_{i=1}^k e_i$$ for any scalars $a_i(e)$ ($i \in [1, k],\;
e \in E'$) and any finite $E' \subseteq E$.
(Thus, $\bigwedge_{i=1}^k e_i = 0$ unless
all $e_i$ are distinct.) 
The inner product on $\Lambda^k E$ satisfies
\rlabel e.ipdet
{
\ip{u_1 \wedge \cdots \wedge u_k, v_1 \wedge \cdots \wedge v_k}
=
\det\big[\ip{u_i, v_j}\big]_{i, j \in [1, k]}
}
when $u_i$ and $v_j$ are 1-vectors.
(This also shows that the inner product on $\Lambda^k E$ does not
depend on the choice of orthonormal basis of $\eH$.)
We then define the \dfnterm{exterior} (or \dfnterm{Grassmann}) \dfnterm{algebra}
$\Ext\big(\eH\big) := \Ext(E) := \bigoplus_{k \ge 0} \Lambda^k E$, where the
summands are declared orthogonal, making it into a Hilbert space.
\iffalse
(Throughout the paper, $\oplus$ is used to indicate the sum
of orthogonal summands, or, if there are an infinite number of orthogonal
summands, the closure of their sum.)
\fi
Vectors $u_1, \ldots, u_k \in \eH$ are linearly independent iff $u_1
\wedge \cdots \wedge u_k \ne 0$.
For a $k$-element subset $A \subseteq E$ with ordering $\Seq{e_i}$ in $E_k$,
write
$
\theta_A := \bigwedge_{i=1}^k e_i
$.
We also write 
$
\bigwedge_{e \in A} f(e) := \bigwedge_{i=1}^k f(e_i)
$
for any function $f \colon E \to \eH$.

Although there is an isometric isomorphism 
$$
u_1 \wedge \cdots \wedge u_k
\mapsto \frac{1}{\sqrt{k!}}
\sum_{\sigma \in \Sym(k)} (-1)^\sigma u_{\sigma(1)} \otimes \cdots
\otimes u_{\sigma(k)} \in \ell^2(E^k)
$$
for $u_i \in \eH$,
this does not simplify matters in the discrete case.
It will be very useful in the continuous case later, however.

If $H$ is a closed linear subspace of $\eH$, written $H \subsp \eH$,
then we identify $\Ext(H)$
with its inclusion in $\Ext(E)$. That is, $\bigwedge^k H$ is the closure of
the linear span of the $k$-vectors
$\{v_1 \wedge \cdots \wedge v_k \st v_1, \ldots, v_k \in
H\}$.
In particular, if $\dim H = r < \infty$, then
$\bigwedge^r H$ is a 1-dimensional subspace of $\Ext(E)$; denote by
$\mv_H$ a unit multivector in this subspace.
Note that $\mv_H$ is unique up to a scalar factor of modulus 1; which
scalar is chosen will not affect the definitions below.
We denote by $P_H$ the orthogonal projection onto $H$ for any 
$H \subsp \eH$ or, more generally, $H \subsp \Ext(E)$. 

\procl l.projection
For every closed 
subspace $H \subsp \eH$, every $k \ge 1$, and every $u_1, \ldots, u_k
\in \eH$,
we have
$
P_{\Ext(H)} (u_1 \wedge \cdots \wedge u_k) = (P_H u_1) \wedge \cdots
\wedge (P_H u_k)
.
$
\endprocl

\iffalse
\rproof 
Write 
$$
u_1 \wedge \cdots \wedge u_k
=
(P_H u_1 + P^\perp_H u_1) \wedge \cdots \wedge (P_H u_k + P^\perp_H u_k)
$$
and expand the product. All terms but $P_H u_1 \wedge \cdots \wedge P_H
u_k$ have a factor of $P^\perp_H u$ in them, making them orthogonal to
$\Ext(H)$ by \rref e.ipdet/. 
\Qed
\fi

\iffalse
A multivector is called \dfnterm{simple} or \dfnterm{decomposable} if it is the
wedge product of 1-vectors.  \rref b.Whitney:book/, p.~49, shows that 
\rlabel e.whitney
{
\|\bu \wedge \bv \| \le \|\bu\| \|\bv\| \quad\hbox{ if either $\bu$ or $\bv$ is
simple}. 
}

We shall use the \dfnterm{interior product} defined by duality:
$$
\ip{\bu \vee \bv, \bw} = \ip{\bu, \bw \wedge \bv}
\qquad (\bu \in \Lambda^{k+l} E ,\ \bv \in \Lambda^l E, \ \bw \in
\Lambda^k E)\,.
$$
In particular, if $e\in E$ and $\bu$ is a multivector that does not contain
any term with $e$ in it (that is, $\bu\in\Ext(e^\perp)$), then $(\bu \wedge
e) \vee e = \bu$ and $\bu \vee e = 0$.
More generally, if $v \in \eH$ with $\|v\|=1$ and $\bu \in \Ext(v^\perp)$,
then $(\bu \wedge v) \vee v = \bu$ and $\bu \vee v = 0$.
Note that the interior product is sesquilinear, not bilinear, over $\bbC$.
\fi

For $v \in \eH$, write $[v]$ for the subspace of scalar multiples of $v$ in
$\eH$.
\iffalse
If $H$ is a finite-dimensional subspace of $\eH$ and $e \notin H$, then 
\rlabel e.Hwedge
{
\mv_H \wedge e = \|P_H^\perp e\|\;\mv_{H + [e]}
}
(up to signum).
To see this, let $u_1, u_2, \ldots, u_r$ be an orthonormal basis of $H$,
where $r = \dim H$.
Put $v := P_H^\perp e /\|P_H^\perp e\|$.
Then $u_1, \ldots, u_r, v$ is an orthonormal basis of $H + [e]$, whence
$$
\mv_{H + [e]} 
=
u_1 \wedge u_2 \wedge \cdots \wedge u_r \wedge v
=
\mv_H \wedge v
=
\mv_H \wedge e/\|P_H^\perp e\|
$$
since $\mv_H \wedge P_H e = 0$.
This shows \rref e.Hwedge/.
Similarly, if $e \notin H^\perp$, then 
\rlabel e.Hvee
{
\mv_H \vee e = \|P_H e\|\;\mv_{H \cap e^\perp}
}
(up to signum).
Indeed, put $w_1 := P_H e /\|P_H e\|$. Let $w_1, w_2, \ldots, w_r$ be an
orthonormal basis of $H$ with $\mv_H = w_1 \wedge w_2 \wedge \cdots \wedge
w_r$.
Then 
$$
\mv_H \vee e
=
\mv_H \vee P_H e
=
(-1)^{r-1} \|P_H e\|\; w_2 \wedge w_3 \cdots \wedge w_r
$$
(up to signum), as desired.

Finally, we claim that 
\rlabel e.reverse
{
\all {u, v \in \eH}\quad \bigip{\mv_H \vee u, \mv_H \vee v} 
=
\ip{P_H v, u}
\,.
}
Indeed, $\mv_H \vee u = \mv_H \vee P_H u$, so this is equivalent to
$$
\bigip{\mv_H \vee P_H u, \mv_H \vee P_H v} 
=
\ip{P_H v, P_H u}
\,.
$$
Thus, it suffices to show that
$$
\all {u, v \in H}\quad \bigip{\mv_H \vee u, \mv_H \vee v} 
=
\ip{v, u}
\,.
$$
By sesquilinearity, it suffices to show this for $u, v$ members of an
orthonormal basis of $H$.
But then it is obvious.

For a more detailed presentation of exterior algebra, see \rref
b.Whitney:book/.
\fi

\bsubsection{Orthogonal Projections}{S.orth}

Let $H$ be a subspace of $\eH$ of dimension $r < \infty$.
Define the probability measure $\bP^H$ on subsets $B \subseteq E$ by
\rlabel e.xiHpr
{
\bP^H\big(\{B\}\big)
:=
|\leftip{\mv_H, \theta_B}|^2
\,.  
}
Note that this is non-0 only for $|B| = r$.
Also, by \rref l.projection/,
$$
\bP^H\big(\{B\}\big) 
=
\|P_{\Ext(H)} \theta_B \|^2 
= 
\| \bigwedge_{e \in B} P_H e \|^2
$$
for $|B| = r$,
which is non-0 iff $\Seq{P_H e \st e\in B}$ are linearly independent.
That is, $\bP^H\big(\{B\}\big) \ne 0$ iff the projections of the elements
of $B$ form a basis of $H$.
%We may also write \rref e.xiHpr/ as
%$$
%\mv_H = \sum_{B \in \cB} \epsilon_B \sqrt{\bP^H[B]} \theta_B
%$$
%for some $\epsilon_B$ of absolute value 1, or alternatively as
%$$
%\bP^H[\ba= B]
%=
%\|P_{\Ext(H)} \theta_{B}\|^2
%=
%\ip{P_{\Ext(H)} \theta_{B}, \theta_{B}}
%\,.
%$$
%
Let $\Seq{v_1, \ldots, v_r}$ be any basis of $H$.
If we use \rref
e.ipdet/ and the fact that $\mv_H = c \bigwedge_i v_i$ for some scalar $c$,
then we obtain 
another formula for $\bP^H$:
\rlabel{e.rowrep}
{
\bP^H\big(\{e_1, \ldots, e_r\}\big) 
  = (\det [\ip{v_{i}, e_{j}}]_{i, j \le r})^2/\det [\ip{v_i, v_j}]_{i, j
  \le r}
\,.  
}

We use $\ba$ to denote a random subset of $E$ arising from a probability
measure $\bP^H$.
To see that \rref e.DPM/ holds for the
matrix of $P_H$, observe that for $|B| = r$,
$$
\bP^H[\ba= B]
=
\bigip{P_{\Ext(H)} \theta_B, \theta_B}
=
\Bigip{\bigwedge_{e \in B} P_H e, \bigwedge_{e \in B} e}
=
\det [\ip{P_H e, f}]_{e, f \in B}
$$
by \rref e.ipdet/.
This shows that \rref e.DPM/ holds for $|A| = r$ since $|\ba| = r$ $\bP^H$-a.s.
The general case is a consequence of multilinearity, which gives
the following extension of \rref e.DPM/.
We use the convention that $\theta_ \emptyset := 1$ and $\bu \wedge 1 :=
\bu$ for any multivector $\bu$.

\procl t.genprs
If $A_1$ and  $A_2$ are (possibly empty) subsets of a finite set $E$, 
then
\rlabel e.genprs
{
\bP^H[A_1 \subseteq \ba, A_2 \cap \ba = \emptyset]
=
\bigip{P_{\Ext(H)} \theta_{A_1} \wedge P_{\Ext(H^\perp)} \theta_{A_2},
\theta_{A_1} \wedge \theta_{A_2}}
\,.
}
In particular, for every $A \subseteq E$, we have
\rlabel e.included
{
\bP^H[A \subseteq \ba]
=
\|P_{\Ext(H)} \theta_A\|^2
\,.
}
\endprocl

%We see immediately the relationship of orthogonality to duality: 

\procl c.dualrep
If $E$ is finite, then for every
subspace $H \subsp \eH$, we have
\rlabel e.dualrep
{
\all {B \subseteq E}\quad \bP^{H^\perp}\big(\{E \setminus B\}\big) =
\bP^H\big(\{B\}\big)
\,.
}
\endprocl

These extend to infinite $E$. In order to define $\bP^H$ when $H$ is
infinite dimensional, we proceed by finite approximation.

Let $E = \{e_i \st i \ge 1\}$ be infinite.
Consider first a finite-dimensional subspace $H$ of
$\eH$. Define $H_k$ as the image of
the orthogonal projection of $H$ onto the span of $\{e_i \st 1 \le i \le k\}$.
By considering a basis of $H$, we see that $P_{H_k} \to P_H$ in the weak
operator topology (WOT), i.e., matrix-entrywise,
as $k \to\infty$. 
It is also easy to
see that if $r := \dim H$, then $\dim H_k = r$ for all large $k$ and, in
fact, $\mv_{H_k} \to \mv_H$ in the usual norm topology. It follows that
\rref e.genprs/ holds for this subspace $H$ and for every finite $A_1,
A_2 \subset E$.

Now let $H$ be an infinite-dimensional closed subspace of $\eH$. 
Choose finite-dimensional subspaces $H_k \uparrow H$.
It is well known that 
$P_{H_k} \to P_H$ (WOT).
Then
\rlabel e.detgenprs
{
\hbox{for all finite sets } A \quad
\det (P_{H_k} \restrict A) \to \det (P_H \restrict A)
\,,
}
whence $\bP^{H_k}$ has a weak${}^*$ limit that we
denote $\bP^H$ and that satisfies
\rref e.genprs/.

We also note that for {\it any\/} sequence of subspaces $H_k$, if $P_{H_k}
\to P_H$ (WOT), then $\bP^{H_k} \to \bP^H$ weak${}^*$ because \rref
e.detgenprs/ then holds.

\bsubsection{Positive Contractions}{S.contract}

We call $Q$ a \dfnterm{positive contraction} if $Q$ is a self-adjoint
operator on $\eH$ such that for all $u \in \eH$, we have $0 \le (Q u, u) \le
(u, u)$.
%To show existence of a corresponding determinantal probability measure, which
%we shall denote $\bP^Q$, let $P_H$ be any orthogonal projection that is
A \dfnterm{projection dilation} of $Q$ is an orthogonal projection $P_H$ onto
a closed subspace $H \subsp \ell^2(E')$ for some $E' \supseteq E$ such
that for all $u \in \eH$, we have $Qu = P_{\eH} P_H u$, where we regard
$\ell^2(E')$ as the orthogonal sum $\eH \oplus \ell^2(E' \setminus E)$. In
this case, $Q$ is also called the \dfnterm{compression} of $P_H$ to $\eH$.
Choose such a dilation
(see \rref e.vecdilate/ or \rref e.dilate/)
%The existence of a dilation is standard and is easily
%constructed (see \rref e.vecdilate/ or \rref e.dilate/).
and define $\bP^Q$ as the law of $\ba \cap
E$ when $\ba$ has the law $\bP^H$.
Then \rref e.DPM/ for $Q$ is a special case of \rref e.DPM/ for $P_H$.

Of course, when $Q$ is the orthogonal projection onto a subspace $H$, then
$\bP^Q = \bP^H$.
Basic properties of $\bP^Q$ follow from those for orthogonal
projections, such as:

\procl t.Q
If\/ $Q$ is a positive contraction, then
for all finite $A_1, A_2 \subseteq E$,
\rlabel e.Qgenprs
{
\bP^Q\left[ A_1 \subseteq \qba, A_2 \cap \qba = \emptyset\right]
=
\Bigip{\bigwedge_{e \in A_1} Q e \wedge \bigwedge_{e \in A_2} (I-Q) e,
\theta_{A_1} \wedge \theta_{A_2} }
\,.
}
\endprocl

If \rref e.DPM/ is given, then \rref e.Qgenprs/ can be deduced from
\rref e.DPM/ without using our general theory and, in fact, without assuming
that the matrix $Q$ is self-adjoint. Indeed,
suppose that $X$ is any diagonal matrix. Denote its $(e, e)$-entry by $x_e$.
Comparing coefficients of $x_e$ shows that \rref e.DPM/ implies, for finite
$A \subseteq E$,
\rlabel e.xe
{
\EBig{\prod_{e \in A} \big(\II{e \in \qba} + x_e\big)}
=
\det \big( (Q + X) \restrict A \big)
\,.
}
Replacing $A$ by $A_1 \cup A_2$ and
choosing $x_e := - \boI{A_2}(e)$ gives \rref e.Qgenprs/.
On the other hand, if
we substitute $x_e := 1/(z_e-1)$, then we may rewrite \rref e.xe/ as
\rlabel e.ze
{
\EBig{\prod_{e \in A} \big(\II{e \in \qba} z_e + \II{e \notin \qba}\big)}
=
\det \big( (Q Z + I-Q) \restrict A \big)
\,,
}
where $Z$ is the diagonal matrix of the variables $z_e$.
Let $E$ be finite.
Write $z^A := \prod_{e \in A} z_e$ for $A \subseteq E$.
Then \rref e.ze/ is equivalent to
\rlabel e.affine
{
\sum_{A \subseteq E} \bP^Q[\qba = A] z^A
=
\det(I-Q+QZ)
\,.
}
This is the same as the Laplace transform of $\bP^Q$ after a trivial change
of variables.
When $\|Q\| < 1$, we can write $\det(I-Q+QZ) = \det(I-Q) \det(I+JZ)$ with
$J := Q(I-Q)^{-1}$. Thus, for all $A \subseteq E$, we have 
\rlabel{e.Jform}
{
\bP^Q[\qba = A]
=
\det(I-Q) \det (J\restrict A)
=
\det(I+J)^{-1} \det (J\restrict A)
\,.
}

A probability measure $\bP$ on $2^E$ is called \dfnterm{strongly Rayleigh}
if its generating polynomial 
$f(z) := \sum_{A \subseteq E} \bP[\qba = A] z^A$ satisfies the inequality
\rlabel{e.strongRay}
{
\frac{\partial f}{\partial z_e}(x)
\frac{\partial f}{\partial z_{e'}}(x)
\ge
\frac{\partial^2 f}{\partial z_e \partial z_{e'}}(x)
f(x)
}
for all $e \ne e' \in E$ and all real $x \in \bbR^E$.
This property is satisfied by every determinantal probability measure, as
was shown by \rref b.BBL:Rayleigh/, who demonstrated its usefulness in
showing other properties, such as negative associations and preservation
under symmetric exclusion processes.

For a set $K
\subseteq E$, denote by $\sF(K)$ the $\sigma$-field of events
that are measurable with respect to the events $\{e \in
\qba\}$ for $e \in K$. 
Define the \dfnterm{tail} $\sigma$-field to be the
intersection of $\sF(E\setminus K)$ over all finite $K$.
We say that a measure $\bP$ on $2^E$ has \dfnterm{trivial tail} if 
every event in the tail $\sigma$-field has measure either 0 or 1.
%Recall that tail triviality is equivalent to
%\rlabel e.georgii
%{
%\forall \cA_1 \in \sF(E) \;\ \forall \epsilon > 0 \;\ \exists K \hbox{
%finite } \;
%\forall \ev A_2 \in \sF(E\setminus K)\qquad
%\bigl|\bP(\cA_1 \cap \ev A_2) - \bP(\cA_1) \bP(\ev A_2)\bigr| < \epsilon\,.
%}
%(See, e.g., \rref b.Georgii:book/, p.~120.)

\rprocl t.tail {\rref b.Lyons:det/}
If\/ $Q$ is a positive contraction, then
$\bP^Q$ has trivial tail.
\endprocl

For finite $E$ and a positive contraction $Q$, define the
\dfnterm{entropy} of $\bP^Q$ to be
$$
\ent(Q) := - \sum_{A \subseteq E} \bP^Q(\{A\}) \log \bP^Q(\{A\}) 
\,.
$$
Numerical calculation supports the following conjecture \rref b.Lyons:det/: 

\procl g.concave
For all positive contractions $Q_1$ and $Q_2$, we have 
\rlabel e.concave
{
\ent\big((Q_1+Q_2)/2\big) \ge \big(\ent(Q_1) + \ent(Q_2)\big)/2
\,.
}
\endprocl

\bsubsection{Stochastic Inequalities}{S.dinequalities}

Let $E$ be denumerable.
A function $f \colon 2^E \to \bbR$ is called \dfnterm{increasing} if for
all $A \in 2^E$ and all $e \in E$, we have $f\big(A \cup \{ e
\}\big) \ge f(A)$.
An event is called increasing or
\dfnterm{upwardly closed} if its indicator is increasing.

Given two probability measures $\bP^1$, $\bP^2$ on $2^E$, we say that
\dfnterm{$\bP^2$ stochastically dominates $\bP^1$} and write $\bP^1 \dom \bP^2$ if
for all increasing events $\ev A$, we have $\bP^1(\ev A) \le \bP^2(\ev A)$.
This is equivalent to $\int f \,d\bP^1 \le \int f \,d\bP^2$ for all bounded
increasing $f$.

A \dfnterm{coupling} of two probability measures $\bP^1$, $\bP^2$ on
$2^E$ is a probability measure $\mu$
on $2^E \times 2^E$ whose coordinate projections are $\bP^1$, $\bP^2$;
it is
\dfnterm{monotone} if
$$
\mu\big\{(\ev A_1, \ev A_2) \st \ev A_1 \subseteq \ev A_2\big\} = 1
\,.$$
By Strassen's theorem \cite{Strassen}, stochastic
domination $\bP^1 \preccurlyeq \bP^2$ is equivalent to the existence of a
monotone coupling of $\bP^1$ and $\bP^2$.

\rprocl t.dominate-infinite {\rref b.Lyons:det/}
If $H_1 \subsp H_2 \subsp \eH$, then $\bP^{H_1} \preccurlyeq
\bP^{H_2}$.
\endprocl

It would be very interesting 
to find a natural or explicit monotone coupling.

%A coupling $\mu$ is \dfnterm{disjoint} 
%if $\mu\big\{(A_1, A_2) \st A_1 \cap A_2 = \emptyset\big\} = 1$.
A coupling $\mu$ has 
\dfnterm{union marginal} $\bP$ if
for all events $\cA \subseteq 2^E$, we have
$
\bP(\cA) = \mu\big\{(A_1, A_2) \st A_1 \cup A_2 \in \cA \big \}
$.

\rprocl q.unioncoupling \cite{Lyons:det}
Given $H = H_1 \oplus H_2$, is there a  
coupling of\/ $\bP^{H_1}$ and $\bP^{H_2}$ with union marginal $\bP^H$?  
\endprocl

A positive answer is supported by some numerical calculation.
It is easily seen to hold when
$H=\eH$
by \rref c.dualrep/.  
%: The probability measure $\mu$ on $2^E \times
%2^E$ defined by 
%$$
%\mu\{(A, E \setminus A) \st A \in \cA\} := \bP^{H_1}(\cA)
%$$
%and 
%$$
%\mu\big \{ (A, B) \st B \ne E \setminus A \} := 0
%$$
%does this, as we can see by \rref c.dualrep/.  

In the sequel, we write $Q_1 \lLoew Q_2$ if $\ip{Q_1 u, u} \le \ip{Q_2 u, u}$
for all $u \in \eH$.
%The following extension of \rref t.dominate-infinite/
%was proved by \rref b.Lyons:det/
%and \rref b.BBL:Rayleigh/:

\rprocl t.dominate \cite{Lyons:det,BBL:Rayleigh}
If\/ $0 \lLoew Q_1 \lLoew Q_2 \lLoew I$, then $\bP^{Q_1} \dom \bP^{Q_2}$.
\endprocl

\rproof
By \rref t.dominate-infinite/, it suffices that
there exist orthogonal projections $P_1$ and $P_2$ that are dilations of
$Q_1$ and $Q_2$ such that $P_1 \lLoew P_2$.
This follows from Na{\u\i}mark's dilation
theorem \cite{Paulsen}, which says that
any measure
whose values are positive operators, whose total mass is $I$, and which is
countably additive in the weak operator topology
dilates to a spectral measure.  The measure in our case is defined on a
3-point space, with masses $Q_1$, $Q_2
-Q_1$, and $I-Q_2$, respectively.
If we denote the respective dilations by $R_1$, $R_2$, and $R_3$, then we
set $P_1 := R_1$ and $P_2 := R_1 + R_2$.
\Qed

A positive answer in general to \rref q.unioncoupling/
would give the following more general
result by compression: If $Q_1$, $Q_2$ and $Q_1 + Q_2$ are
positive contractions on $\eH$, then
there is a coupling of $\bP^{Q_1}$ and $\bP^{Q_2}$
with union marginal $\bP^{Q_1+Q_2}$.
%We note that the requirement of being disjoint is superfluous.

It would be very useful to have additional sufficient conditions for
stochastic domination:
see the end of \rref S.orthogpoly/ and \rref g.FKdom/.
For examples where more is known, see \rref t.GMdom/. 

We shall say that the events in $\sF(K)$ are \dfnterm{measurable with
respect to} $K$ and likewise for functions that are measurable with respect
to $\sF(K)$.
We say that $\bP$ has \dfnterm{negative associations} if
for every pair $f_1$, $f_2$ of increasing functions that are measurable with
respect to complementary subsets of $E$,
\rlabel e.negass
{
\bE[f_1 f_2] \le \bE[f_1] \bE[f_2] 
\,.
}
%In this case, for any collection $f_1, f_2, \ldots, f_n$ of increasing
%{\it nonnegative\/}
%functions that are measurable with respect to pairwise disjoint
%subsets of $E$, we have
%\rlabel e.negass-multi
%{
%\bE[f_1 f_2 \cdots f_n] \le \bE[f_1] \bE[f_2] \cdots \bE[f_n]
%\,.
%}
%This is shown by an easy induction argument.
%One could replace ``increasing" by ``decreasing" just as well.

%Extending a result of
%\cite[Lemma 3.2]{FedMih},
%\rref b.Lyons:det/ proved the following (for the meaning of the last terms,
%see \rref b.Lyons:det/):

\rprocl t.FM \cite{Lyons:det}
If\/ $0 \lLoew Q \lLoew I$, then
$\bP^Q$ has negative associations.
%and, in fact, conditional negative
%associations with external fields.
\endprocl

\rproof
The details for finite $E$ were given in \rref b.Lyons:det/.
For infinite $E$, let
$f_1$ and $f_2$ be increasing bounded functions
measurable with respect to
$\sF(A)$ and $\sF(E \setminus A)$,
respectively. Choose finite $E_n
\uparrow E$. The
conditional expectations $\bE[f_1 \mid \sF(A \cap E_n)]$ and
$\bE[f_2 \mid \sF(E_n \setminus A)]$ are increasing functions to which
\rref e.negass/ applies (because restriction to $E_n$ corresponds to a
compression of $Q$, which is a positive contraction) and which, being
martingales,
converge to
$f_1$ and $f_2$ in $L^2(\bP^Q)$.
\Qed

\bsubsection{Mixtures}{S.mix}

%We call a self-adjoint operator $T$ on $\eH$ \dfnterm{trace class} if it is
%compact and its eigenvalues are absolutely summable.
Write $\Bern(p)$ for the distribution of a Bernoulli random variable with
expectation $p$.
For $p_k \in [0, 1]$, let $\Bin(\Seq{p_k})$ be the distribution of a sum
of independent $\Bern(p_k)$ random variables.
Recall that $[v]$ is the set of  scalar multiples of $v$.

%The following fact was first observed by \rref b.Bapat/ for the case of
%uniform spanning trees.

%\procl t.eigcount
%Let $Q$ be a positive contraction with eigenvalues
%$\Seq{\lambda_k \st k \ge 1}$, listed with multiplicity.
%If $\qba \sim\bP^Q$, then the distribution of $|\qba|$ is
%$\Bin(\Seq{\lambda_k})$.
%\endprocl
%
%This is an immediate consequence of the following
%extension, which was first observed by \cite[Lemma 3.4]{ShiTak:I} and
%\cite[(2.38)]{ShiTak:II} in a less probabilistic form.
%It was made a central part of the construction of determinantal probability
%measures corresponding to positive contractions 
%by \rref b.HKPV:survey/, who showed how useful it is.

\rprocl t.eigmix {\cite{Bapat}; Lemma 3.4 of \cite{ShiTak:I}; (2.38) of
\cite{ShiTak:II}; \cite{HKPV:survey}}
Let $Q$ be a positive contraction with spectral decomposition
$Q = \sum_k \lambda_k P_{[v_k]}$,
where $\Seq{v_k \st k \ge 1}$ are orthonormal.
Let $I_k \sim \Bern(\lambda_k)$ be independent.
Let $\rH := \bigoplus_k [I_k v_k]$; thus, $Q = \bE P_{\rH}$.
Then $\bP^Q = \bE \bP^{\rH}$.
Hence, if\/ $\qba \sim\bP^Q$, then $|\qba| \sim \Bin(\Seq{\lambda_k})$.
\endprocl

\rproof
%There are various ways to prove this.
By \rref t.dominate/, 
it suffices to prove it when only finitely many $\lambda_k \ne 0$.
Then by \rref t.Q/, we have
$
\bP^Q\left[ A \subseteq \qba \right]
=
\Bigip{\bigwedge_{e \in A} Q e, \theta_A }
$
for all $A \subseteq E$.
Now
\begin{align*}
\bigwedge_{e \in A} Q e
=
\bigwedge_{e \in A} \sum_k \lambda_k P_{[v_k]} e 
&=
\sum_{j \colon A \to \bbN} \prod_{e \in A} \lambda_{j(e)} 
\bigwedge_{e \in A} P_{[v_{j(e)}]} e
\\
&=
\sum_{j \colon A \rightarrowtail \bbN} \prod_{e \in A} \lambda_{j(e)}
\bigwedge_{e \in A} P_{[v_{j(e)}]} e
\end{align*}
because $v \wedge v = 0$ and $P_{[v]} e$ is a multiple of $v$, so none of
the terms where $j$ is not injective contribute.
Thus,
\begin{align*}
\bigwedge_{e \in A} Q e
&=
\sum_{j \colon A \rightarrowtail \bbN} \EBig{\prod_{e \in A} I_{j(e)}}
\bigwedge_{e \in A} P_{[v_{j(e)}]} e
=
\EBig{\sum_{j \colon A \rightarrowtail \bbN} \prod_{e \in A} I_{j(e)}
\bigwedge_{e \in A} P_{[v_{j(e)}]} e}
\\
&=
\EBig{\sum_{j \colon A \to \bbN} \prod_{e \in A} I_{j(e)}
\bigwedge_{e \in A} P_{[v_{j(e)}]} e}
=
\bE\bigwedge_{e \in A} \sum_k I_k P_{[v_k]} e 
=
\bE \bigwedge_{e \in A} P_{\rH} e 
\,.
\end{align*}
We conclude that
$\bP^Q\left[ A \subseteq \qba \right]
=
\bE \leftip{\bigwedge_{e \in A} P_{\rH} e , \theta_A}
=
\Ebig{ \bP^{\rH}\left[ A  \subseteq \ba \right]}
$
by \rref e.Qgenprs/. 
\Qed

We sketch another proof:
Let $E'$ be disjoint from $E$ with the same cardinality. 
Choose an orthonormal sequence $\Seq{v_k'}$ in $\ell^2(E')$.
Define 
\rlabel{e.vecdilate}
{
H := \bigoplus_k \big[\sqrt{\lambda_k} v_k + \sqrt{1 - \lambda_k} v_k'\big]
\subsp \ell^2(E \cup E')
\,.
}
Then $Q$ is the compression of $P_H$ to $\eH$.
Expanding $\mv_H = \bigwedge_k (\sqrt{\lambda_k} v_k + \sqrt{1 - \lambda_k}
v_k')$ in the obvious way into orthogonal pieces and restricting to $E$, we
obtain the desired equation from \rref e.xiHpr/.

The first proof shows more generally the following:
Let $Q_0$ be a positive contraction.
Let $\Seq{v_k \st k \ge 1}$ be (not necessarily orthogonal) vectors
such that $Q_0 + \sum_k P_{[v_k]} \lLoew I$.
Let $I_k$ be independent Bernoulli random variables
with $\bE \sum_k I_k < \infty$.
Write $\rQ := Q_0 + \sum_k I_k P_{[v_k]}$.
Then $\bP^{\bE \rQ} = \bE \bP^\rQ$.
This was observed by Ghosh and Krishnapur (personal communication, 2014).

%One may ask about other mixtures of determinantal probability measures that
%correspond to orthogonal projections.
Note that in the mixture of \rref t.eigmix/, the distribution of $\Seq{I_k
\st k \ge 1}$ is determinantal corresponding to the diagonal matrix with
diagonal $\Seq{\lambda_k \st k \ge 1}$.
Thus, it is natural to wonder whether $\Seq{I_k \st k \ge 1}$ can be taken
to be a general determinantal measure. If such a mixture is not necessarily
determinantal, must it be strongly Rayleigh or at least have negative
correlations?
Here, we say that a probability measure $\bP$ on $2^E$ has
\dfnterm{negative correlations} if for every pair $A$, $B$ of finite
disjoint subsets of $E$, we have 
$
\bP[ A \cup B \subseteq \qba] 
\le
\bP[ A \subseteq \qba] \bP[ B \subseteq \qba] 
$.
Note that negative associations is stronger than negative correlations.

\bsubsection{Example: Uniform Spanning Trees and Forests}{S.ust}

The most well-known example of a (nontrivial discrete) determinantal
probability measure is that where $\qba$ is a uniformly chosen random spanning
tree of a finite connected graph $G = (\vertex, \edge)$ with $E := \edge$.
Here, we regard a spanning tree as a set of edges.
%In this case, $Q$ is the \dfnterm{transfer current matrix} $Y$,
%which is defined as follows.
%Orient the edges of $G$ arbitrarily.
%Regard $G$ as an electrical network with each edge having unit
%conductance.  Then $Y(e, f)$ is the amount of current flowing along the edge
%$f$ when a battery is hooked up between the endpoints of $e$ of such
%voltage that in the network as a whole, unit current flows from the tail of
%$e$ to the head of $e$.
The fact that \eqref{e.DPM} holds for the uniform spanning tree is due to
\rref b.BurPem/ and is called the Transfer Current Theorem.
The case with $|A| = 1$ was shown much earlier by \rref b.Kirchhoff/, while
the case with $|A| = 2$ was first shown by \rref b.BSST/.
Write $\ust_G$ for the uniform spanning tree measure on $G$.

To see that $\ust_G$ is indeed determinantal, consider the
vertex-edge incidence matrix $M$ of $G$, where each edge is oriented
(arbitrarily) and the $(x, e)$-entry of $M$ equals 1 if $x$ is the head of
$e$, $-1$ if $x$ is the tail of $e$, and 0 otherwise.
Identifying an edge with its corresponding column of $M$,
we find that a spanning tree is
the same as a basis of the column space of $M$. Given $x \in \vertex$,
define the \dfnterm{star} at $x$ to be the $x$-row of $M$, regarded as a
vector $\str_x$ in the row space, $\STAR(G) \subsp \ell^2(\edges)$.
It is easy that the row-rank of $M$ is $|\verts| - 1$.
Let $x_0 \in \verts$ and let $\bu$ be the
wedge product (in some order) of the stars at all the vertices other than
$x_0$. Thus, $\bu = c\, \mv_{\STAR(G)}$ for some $c \ne 0$.
Since spanning trees are bases of the column space of $M$, we have
$\bigip{\bu, \theta_A} \ne 0$ iff $A$ is a spanning tree.
That is, the only non-zero coefficients of $\bu$ are those in
which choosing one edge in each $\str_x$ for $x \ne x_0$ yields a spanning
tree; moreover, each spanning tree occurs exactly once since there is exactly
one way to choose an edge incident to each $x \ne x_0$ to get a given spanning
tree. 
This means that its coefficient is $\pm 1$. 
Hence, $\bP^{\STAR(G)}$ is indeed uniform on spanning trees.
Simultaneously, this proves the matrix tree theorem that the
number of spanning trees equals $\det[\ip{\str_x, \str_y}]_{x, y \ne x_0}$,
since this determinant is $\|\bu\|^2$.

One can define analogues of $\ust_G$ on infinite connected graphs
\cite{Pemantle:ust,Hag:rcust,BLPS:usf} by weak limits.
For brevity, we simply define them here as determinantal probability
measures.
Again, all edges of $G$ are oriented arbitrarily.
We define $\STAR(G)$ as the closure of the linear span of the stars. An
element of $\ell^2(\edges)$ that is finitely supported and orthogonal to
$\STAR(G)$ is called a \dfnterm{cycle}; the closed linear span of the
cycles is $\CYCLE(G)$.
The \dfnterm{wired uniform spanning forest} is
$\wsf_G := \bP^{\STAR(G)}$, 
while 
the \dfnterm{free uniform spanning forest} 
is $\fsf_G := \bP^{\CYCLE(G)^\perp}$.
%In particular, $\wsf_G = \fsf_G$ iff $\STAR(G) = \CYCLE(G)^\perp$.

\bsection{Continuous Basics}{s.continuous}

Our discussion of the ``continuous" case includes the discrete case, but
the discrete case has the more elementary formulations given earlier.

Let $E$ be a measurable space. As before, $E$ will play the role of the
underlying set on which a point process forms a counting measure.
While before we implicitly used counting measure on $E$ itself, now we
shall have an arbitrary measure $\mu$; it need not be a probability
measure. The case of Lebesgue measure on Euclidean space is a common
one.
The Hilbert spaces of interest will be $L^2(E, \mu)$.

\bsubsection{Symmetrization and Anti-symmetrization}{S.symm}

There may be no natural order in $E$, so to define, e.g., a probability
measure on $n$ points of $E$, it is natural to use a probability measure on
$E^n$ that is symmetric under coordinate changes and that vanishes on the
diagonal $\diag_n(E) := \big\{(x_1, \ldots, x_n) \in E^n \st \texists {i \ne j}
x_i = x_j\big\}$.
Likewise, for exterior algebra, it is more convenient to identify $u_1
\wedge \cdots \wedge u_n$ with $$\frac{1}{\sqrt{n!}}\sum_{\sigma \in \Sym(n)}
(-1)^\sigma u_{\sigma(1)}
\otimes \cdots \otimes u_{\sigma(n)} \in L^2(E^n, \mu^n)$$ for
$u_i \in L^2(E, \mu)$.
Thus, $u_1 \wedge \cdots \wedge u_n$ is identified with the function 
$
(x_1, \ldots, x_n) \mapsto 
%\frac{1}{\sqrt{n!}}
\det [u_i(x_j)]_{i, j \in \{1, \ldots, n\}}/\sqrt{n!}
$.
Note that
\begin{align} \label{e.prodtensor}
n! \Big(\bigwedge_{i=1}^n u_i\Big) \Big(\bigwedge_{i=1}^n v_i\Big)(x_1,
\ldots, x_n)
&=
\det [u_i(x_j)] %_{i, j \in \{1, \ldots, n\}}
\det [{v_i(x_j)}] %_{i, j \in \{1, \ldots, n\}}
=
\det [u_i(x_j)] %_{i, j \in \{1, \ldots, n\}}
\det [{v_i(x_j)}]^T %_{i, j \in \{1, \ldots, n\}}
\nonumber \\
&=
\det [u_i(x_j)][{v_i(x_j)}]^T
=
\det [K(x_i, x_j)]_{i, j \in \{1, \ldots, n\}}
\end{align}
with $K := \sum_{i=1}^n u_i \otimes {v_i}$.
Here, ${}^T$ denotes transpose.

\bsubsection{Joint Intensities}{S.intense}

Suppose from now on that $E$ is a locally compact Polish space
(equivalently, a locally compact second countable Hausdorff space).
Let $\mu$ be a Radon measure on $E$, i.e., a Borel measure that is finite
on compact sets.
Let $\cN(E)$ be the set of Radon measures on $E$ with values in $\bbN \cup
\{\infty\}$.
We give $\cN(E)$ the vague topology generated by the maps $\xi \mapsto
\int f \,d\xi$ for continuous $f$ with compact support; then $\cN(E)$ is
Polish.
The corresponding Borel $\sigma$-field of
$\cN(E)$ is generated by the maps $\xi \mapsto
\xi(A)$ for Borel $A \subseteq E$.
Let $\sX$ be a simple point process on $E$, i.e., a random variable with
values in $\cN(E)$ such that $\sX(\{x\}) \in \{0, 1\}$ for all $x \in E$.
The power $\sX^k := \sX \otimes \cdots \otimes \sX$ lies in $\cN(E^k)$.
%We always assume that $\bE[\sX]$ is a Radon measure.
Thus, $\bE[\sX^k]$ is a Borel measure on $E^k$; the part of it that is
concentrated on $E^k \setminus \diag_k(E)$ is called the \dfnterm{$k$-point
intensity measure} of $\sX$. If the intensity measure is absolutely
continuous with respect to $\mu^k$, then its Radon-Nikodym derivative
$\rho_k$ is called the \dfnterm{$k$-point intensity function} or the
\dfnterm{$k$-point correlation function}:
\rlabel{e.RN}
{
\mbox{for all Borel } A \subseteq E^k \setminus \diag_k(E)
\qquad
\bE[\sX^k(A)]
=
\int_A \rho_k \,d\mu^k
\,.
}
Since the intensity measure vanishes on the diagonal $\diag_k(E)$, we take
$\rho_k$ to vanish on $\diag_k(E)$.
We also take $\rho_k$ to be symmetric under permutations of coordinates.
Intensity functions are the continuous analogue of the elementary
probabilities \rref e.DPM/.

Since the sets $\prod_{i=1}^k A_i := A_1 \times \cdots \times A_k$ generate
the $\sigma$-field on $E^k \setminus \diag_k(E)$ for pairwise
disjoint Borel $A_1, \ldots, A_k \subseteq E$, a measurable function $\rho_k
\colon E^k \to [0, \infty)$ is ``the" $k$-point intensity function iff 
\rlabel{e.disj}
{
\EBig{\prod_{i=1}^k \sX(A_i)}
=
\int_{\prod_{i=1}^k A_i} \rho_k \,d\mu^k
\,.
}
Since $\sX$ is simple,
$\sX^k\big(A^k \setminus D_k(A)\big) = \big(\sX(A)\big)_k$,
where $(n)_k := n(n-1) \cdots (n-k+1)$.
%% This uses that $\sX$ is simple.
Since $\rho_k$ vanishes on the diagonal, it follows from \rref e.RN/
that for disjoint $A_1,
\ldots, A_r$ and non-negative $k_1, \ldots, k_r$ summing to $k$,
\rlabel{e.falling}
{
\EBig{\prod_{j=1}^r \big(\sX(A_j)\big)_{k_j}}
=
\int_{\prod_{j=1}^r A_j^{k_j}} \rho_k \,d\mu^k
\,.
}
Again, this characterizes $\rho_k$, even if we use only $r = 1$.

In the special case that $\sX(E) = n$ a.s.\ for some $n \in \bbZ^+$, then
the definition \rref e.RN/ shows that a random ordering of the $n$ points
of $\sX$ has density $\rho_n/n!$. More generally, \rref e.RN/ shows that
for all $k < n$, 
\rlabel{e.densityktuple}
{
\mbox{the density of a random (ordered) $k$-tuple of $\sX$
is $\rho_k/(n)_k\,,$}
}
whence in this case,
\rlabel{e.integrate}
{
\rho_k(x_1, \ldots, x_k)
=
\frac{1}{(n-k)!} \int_{E^{n-k}} \rho_n(x_1, \ldots, x_n)
\,d\mu^{n-k}(x_{k+1}, \ldots, x_n)
\,.
}

We call
$\sX$ \dfnterm{determinantal} if for some measurable $K
\colon E^2 \to \bbC$ and all $k \ge 1$, $\rho_k(F) = \det (K \restrict F)$
$\mu^k$-a.e.
Here, $K \restrict (x_1, \ldots, x_k)$ is the matrix $[K(x_i, x_j)]_{i, j
\le k}$.
In this case, we denote the law of $\sX$ by $\bP^K$.

We consider only $K$ that are locally square integrable (i.e., $|K|^2
\mu^2$ is Radon), are Hermitian (i.e., $K(y, x) = \overline{K(x, y)}$ for
all $x, y \in E$), and are positive semidefinite (i.e., $K \restrict F$ is
positive semidefinite for all finite $F$, written $K \gLoew 0$).
In this case, $K$ defines a positive semidefinite integral operator 
$
(Kf)(x) := \int K(x, y) f(y) \,d\mu(y)
$
on functions $f \in L^2(\mu)$ with compact support.
For every Borel $A \subseteq E$, we denote by $\mu_A$ the measure $\mu$
restricted to Borel subsets of $A$ and by
$K_A$ the compression of $K$ to $A$, i.e., $K_A f := (K f)
\restrict A$ for $f \in L^2(A, \mu_A)$.
The operator $K$ is locally trace-class, i.e.,
for every compact $A \subseteq E$,
the compression $K_A$ 
is trace class, having a
spectral decomposition 
$\medstrut
K_A f = \sum_k \lambda^A_k \ip{f, \phi^A_k} \phi^A_k
$,
where $\Seq{\phi^A_k \st k \ge 1}$ are orthonormal eigenfunctions of $K_A$
with positive summable eigenvalues $\Seq{\lambda^A_k \st k \ge 1}$.
If $A_1$ is the set where $\sum_k \lambda^A_k |\phi^A_k|^2 < \infty$, then
$\mu(A \setminus A_1) = 0$ and
$\medstrut\sum_k \lambda^A_k \phi^A_k \otimes \overline{\phi^A_k}$ converges
on $A_1^2$, with sum $\mu_A^2$-a.e.\ equal to $K$.
We normally redefine $K$ on a set of measure 0 to equal this sum.
Such a $K$ defines a determinantal point process iff the integral operator $K$
extends to all of $L^2(\mu)$ as
a positive contraction \cite{Macchi, Soshnikov:survey, HKPV:survey}.
% Take an increasing sequence of compact subsets; restricted to each we
% have a positive contraction. This means $K$ maps $L^2(\mu)$ to itself
% with norm at most 1.
The joint intensities determine uniquely the law of the point process
\cite[Lemma 4.2.6]{HKPV:book}.
Poisson processes are not determinantal processes,
but when $\mu$ is continuous, they are distributional limits of determinantal
processes.

\bsubsection{Construction}{S.construct}

To see that a positive contraction defines a determinantal point
process, we first consider $K$ that defines an orthogonal projection
onto a finite-dimensional subspace, $H$.
Then $K = \sum_{k=1}^n \phi_k \otimes \overline{\phi_k}$ for every
orthonormal basis $\Seq{\phi_k \st k \le n}$ of $H$ and
$\mv_H = \bigwedge_{i=1}^n \phi_k$ is a
unit multivector in the notation of \rref S.ext/.
Because of \rref e.prodtensor/, we have 
\rlabel{e.firstdensity}
{
\frac{1}{n!} \int \det [K(x_i, x_j)]_{i, j \le n} \,d\mu^n(x_1, \ldots, x_n) 
=
\Big\|\bigwedge_{k=1}^n \phi_k\Big\|^2
=
1
\,,
}
i.e., $\det [K(x_i, x_j)]/n!$ is a density with respect to $\mu^n$.
Although in the discrete case, the absolute squared coefficients of
$\bigwedge_{k=1}^n \phi_k$ give the elementary probabilities, now
coefficients are replaced by a function whose absolute square gives
a probability density.
As noted already, \rref e.firstdensity/
means that $F \mapsto \det (K \restrict F)$ is the
$n$-point intensity function.
In order to show that this density gives a determinantal process with
kernel $K$, we use the
Cauchy-Binet formula, which may be stated as follows:
For $k\times n$ matrices $a = [a_{i, j}]$ and $b = [b_{i, j}]$ with $a^J :=
[a_{i, j}]_{\substack{i \le k \\ j \in J}}$, we have
$$
\det \big([a_{i, j}] [b_{i, j}]^T\big)
=
\sum_{|J| = k} \det a^J \cdot \det b^J
=
\sum_ {\substack{\sigma, \tau \in \Sym(k, n) \\ \img (\sigma) = \img (\tau)}}
(-1)^\sigma (-1)^\tau
\prod_{i=1}^k a_{i, \sigma(i)} b_{i, \tau(i)}
\,,
$$
where $\img(\sigma)$ denotes the image of $\sigma$
and the sums extend over all pairs of injections $$\sigma, \tau \colon
\{1, 2, \ldots, k\} \rightarrowtail \{1, 2, \ldots, n\}\,.$$
Here, the sign $(-1)^\sigma$ of $\sigma$ is defined in the usual way by
the parity of the number of pairs $i < j$ for which $\sigma(i) >
\sigma(j)$.
We have
\begin{align}
\label{e.expanddet}
\rho_k(x_1, \ldots, x_k)
&=
\frac{1}{(n-k)!} \int_{E^{n-k}} \det [K(x_i, x_j)]
\,d\mu^{n-k}(x_{k+1}, \ldots, x_n)
\nonumber \\
&=
\frac{1}{(n-k)!} \int_{E^{n-k}} 
\sum_{\sigma \in \Sym(n)} (-1)^\sigma \prod_{i=1}^n \phi_{\sigma(i)}(x_i)
\cdot {}
\nonumber \\
\noalign{$\displaystyle \hfill \cdot
\sum_{\tau \in \Sym(n)} (-1)^\tau \prod_{i=1}^n \overline{\phi_{\tau(i)}(x_i)}
\,d\mu^{n-k}(x_{k+1}, \ldots, x_n)
$}
&=
\sum_{\substack{\sigma, \tau \in \Sym(k, n) \\ \img (\sigma) = \img (\tau)}}
(-1)^\sigma (-1)^\tau \prod_{i=1}^k \phi_{\sigma(i)}(x_i)
\overline{\phi_{\tau(i)}(x_i)}
\nonumber \\
&=
\det \big(K \restrict (x_1, \ldots, x_k)\big)
\,.
\end{align}
Here, the first equality uses \rref e.integrate/, the second equality uses
\rref e.prodtensor/, the third equality uses the fact that $\int_E
\phi_{\sigma(i)}(x_i) \overline{\phi_{\tau(i)}(x_i)} \,d\mu(x_i)$ is 1 or 0
according as $\sigma(i) = \tau(i)$ or not, and the fourth equality uses
Cauchy-Binet.
Note that a factor of $(n-k)!$ arises because for every pair of injections
$\sigma_1, \tau_1 \in \Sym(k, n)$
with equal image, there are $(n-k)!$ extensions of them to permutations
$\sigma, \tau \in \Sym(n)$
with $\sigma(i) = \tau(i)$ for all $i > k$; in this case,
$(-1)^\sigma (-1)^\tau = (-1)^{\sigma_1} (-1)^{\tau_1}$.
We write $\bP^H$ for the law of the associated point process on $E$.

\procl l.weaklimit
Let $\sX_n \sim \bP^{K_n}$ with $K_n(x, x) \le f(x)$ for some $f \in
L^1_{\rm loc}(E, \mu)$.
Then $\{\bP^{K_n} \st n \ge 1\}$ is tight
and every weak limit point of $\sX_n$ is simple.
%point process.
\endprocl

\rproof
By using the kernel $K_n(x, y)/\sqrt{f(x) f(y)}$ with respect to the
measure $f \mu$, we may assume that $f
\equiv 1$.
Tightness follows from $$m \bP[\sX_n(A) \ge m] \le \bE[\sX_n(A)] 
= \int_A K_n(x, x) \,d\mu(x)\,.$$
For the rest, we may assume that $E$ is compact and $\mu(E) = 1$.
Let $\sX$ be a limit point of $\sX_n$.
Let $\mud$ be the atomic part of $\mu$ and $\muc := \mu - \mud$.
Choose $m \ge 1$ and partition $E$ into sets $A_1, \ldots, A_m$ with
$\muc(A_i) \le 1/m$.
Let $A$ be such that $\mud(E \setminus A) = 0$ and $\muc(A) = 0$. Let $U$
be open such that $A \subseteq U$ and $\muc(U) < 1/m$.
Then 
\begin{align*}
\bP[\sX \mbox{ is not simple}] 
&\le \limsup_n \big(\bP[\sX_n(U \setminus
A) \ge 1] + \bP[\texists i \sX_n(A_i) \ge 2]\big) 
\\ &\le
\limsup_n
\big(\bE[\sX_n(U \setminus A)] + \sum_i \bE[(\sX_n(A_i))_2]\big) 
\\ &\le
\muc(U) + \sum_i \mu(A_i)^2 < 2/m
\,.
\tag*{\qedhere}
\end{align*}
\Qed

Now, given any locally trace-class orthogonal projection $K$ onto
$H$, choose %a sequence of 
finite-dim\-en\-sion\-al subspaces $H_n \uparrow H$ with
corresponding projections $K_n$.
Clearly $K_n(x, y) \to K(x, y)$ $\mu^2$-a.e.\
and $K_n(x, x) \le K(x, x)$ $\mu$-a.e.
Thus, the joint intensity functions converge a.e.
By dominated convergence, 
if $A \subset E^k \setminus \diag_k(E)$ is relatively compact and Borel, then
$\bE^{H_n}[\sX(A)] \to \int_A \det (K \restrict F) \,d\mu^k(F)$.
By uniform exponential moments of $\sX(A)$ \cite[proof of Lemma
4.2.6]{HKPV:book}, it follows that all weak limit points of $\bP^{H_n}$ are
equal, and hence, by \rref l.weaklimit/,
define $\bP^H$ with kernel $K$.
(In \rref S.cinequalities/, we shall see that $\Seq{\bP^{H_n} \st n \ge 1}$
is stochastically increasing.)

Finally, let $K$ be any locally trace-class positive contraction.
Define the orthogonal projection on $L^2(E, \mu) \oplus L^2(E, \mu)$
whose block matrix is 
\rlabel{e.dilate}
{
\begin{pmatrix}
K& \sqrt{K (I-K)} \\  
\sqrt{K (I-K)} &I-K
\end{pmatrix}
.
}
Take an isometric isomorphism of $L^2(E, \mu)$ to $\ell^2(E')$ for some
denumerable set $E'$ and interpret the above as an orthogonal projection
$K'$ on $L^2(E, \mu) \oplus \ell^2(E')$. 
Then $K'$ is clearly locally trace-class and $K$ is the
compression of $K'$ to $E$.
Thus, we define $\bP^K$ by intersecting samples of $\bP^{K'}$ with $E$.
We remark that by
writing $K'$ as a
limit of increasing finite-rank projections that we then compress, we see
that $\bP^K$ may be defined as a limit of determinantal processes
corresponding to increasing finite-rank positive contractions.

\procl g.ctail
If $K$ is a locally trace-class positive contraction, then
$\bP^K$ has trivial tail in that every event in $\bigcap_{{\rm compact\ } A
\subset E} \sF(E \setminus A)$ is trivial.
\endprocl

\bsubsection{Mixtures}{S.cmix}

Rather than using compressions as in the last paragraph above,
an alternative approach to defining $\bP^K$ 
uses mixtures and starts from finite-rank projections, as in \rref S.mix/.
This approach is due to \rref b.HKPV:survey/.
Consider first a finite-rank $K := \sum_{j = 1}^n \lambda_j \phi_j \otimes
\overline{\phi_j}$. 
Let $I_j \sim \Bern(\lambda_j)$ be independent.
Let $\rH := \bigoplus_j [I_j \phi_j]$; thus, $K = \bE P_{\rH}$.
We claim that $\bP^K := \bE \bP^{\rH}$ is determinantal with kernel $K$.
Indeed, it is clearly a simple point process.
Write $\Phi_J := \bigwedge_{j \in J} I_j \phi_j$,
$\psi_j := \sqrt{\lambda_j} \phi_j$, and 
$\psi_J := \bigwedge_{j \in J} \psi_j$.
Let $F \in E^k$.
Combining Cauchy-Binet with \rref e.prodtensor/ yields 
$
\det \big(K \restrict F\big) 
=
k! \sum_{|J| = k} 
|\psi_J(F)|^2 
$.
Similarly,
the joint intensities of $\bE \bP^{\rH}$ are the expectations of
the joint intensities of $\bP^{\rH}$, which equal
\[
\Ebig{\det (P_{\rH} \restrict F)}
=
\EBig{k! \sum_{|J| = k} |\Phi_J(F)|^2 }
=
\det \big(K \restrict F\big) 
\,.
\]

Essentially the same works for trace-class
$K = \sum_{j=1}^\infty \lambda_j \phi_j \otimes \overline{\phi_j}$; we need
merely
take, in the last step, a limit in the above equation as $n \to\infty$ for
$K_n := \sum_{j=1}^n \lambda_j \phi_j \otimes
\overline{\phi_j}$, since all terms are non-negative and $K_n \to K$ a.e.
%In fact, we see that again $\bP^K = \bE \bP^{\rH}$, where $\rH$ is
%finite-dimensional a.s\ by the Borel-Cantelli lemma.
%In \rref S.cinequalities/, we shall see that the processes corresponding to the
%partial sums are stochastically increasing.

Given this construction of $\bP^K$ for trace-class $K$, one can then
construct $\bP^K$ for
a general locally trace-class positive contraction by defining its
restriction to each relatively compact set $A$ via the trace-class
compression $K_A$.

As noted by \rref b.HKPV:survey/,
a consequence of the mixture representation 
is a CLT due originally to \rref b.Soshnikov:Gauss/:

\procl t.CLT
Let $K_n$ be trace-class positive contractions on spaces $L^2(E_n, \mu_n)$.
Let $\sX_n \sim \bP^{K_n}$ and write $|\sX_n| := \sX_n(E_n)$.
If $\Var\big(|\sX_n|\big) \to\infty$ as $n \to\infty$, then $\Seq{|\sX_n|
\st n \ge 1}$ obeys a CLT.
\endprocl

\bsubsection{Simulation}{S.simulate}

In order to simulate $\bP^K$ when $K$ is a trace-class positive
contraction, it suffices, by taking a mixture as above, to see how to
simulate $\sX \sim \bP^H$ when $n := \dim H < \infty$. The following algorithm
\cite[Algo.~18]{HKPV:survey} gives a uniform random ordering of $\sX$ as
$\Seq{X_1, \ldots, X_n}$.
Since $\Ebig{\sX(E)} = n$, the measure $\bE[\sX]/n = n^{-1} K(x,
x)\,d\mu(x)$ is a probability
measure on $E$. Select a point $X_1$ at random from that measure. If $n = 1$,
then we are done. If not, then let $H_1$ be the orthogonal complement in $H$
of the function $K_{X_1} := \sum_{k=1}^n \overline{\phi_k(X_1)}
\phi_k \in H$, where $\Seq{\phi_k \st k \le n}$ is an orthonormal basis for
$H$. Then $\dim H_1 = n-1$ and we may repeat the above for $H_1$ to get the
next point, $X_2$, then $H_2 := H_1 \cap K_{X_2}^\perp$,
etc. The conditional density of $X_{k+1}$
given $X_1, \ldots, X_k$ is $(n-k)^{-1} \det \big(K\restrict (x, X_1, \ldots,
X_k)\big)/ \det \big(K\restrict (X_1, \ldots, X_k)\big)$ by
\rref e.densityktuple/, i.e., $(n-k)^{-1}$ times the squared distance from
$K_x$ to the linear span of $K_{X_1}, \ldots, K_{X_k}$.
It can help for rejection sampling to note that this is at most $(n-k)^{-1}
K(x, x)$. One can also sample faster by noting that the conditional
distribution of $X_{k+1}$ is the same as that of $\bP^\rv$, where $\rv$ is a
uniformly random vector on the unit sphere of $H_k$.

\bsubsection{Transference Principle}{S.transf}

Note that if $N_1, \ldots, N_r$ are bounded $\bbN$-valued random variables,
then the function
$(k_1, \ldots, k_r) \mapsto
\EBig{\prod_{j=1}^r \big(N_j\big)_{k_j}}$ determines the joint distribution
of $\Seq{N_j \st j \le r}$ since it gives the derivatives at $(1, 1,
\ldots, 1)$ of the probability generating function $(s_1, \ldots, s_r)
\mapsto \EBig{\prod_{j=1}^r s_j^{N_j}}$.

Let us re-examine \rref e.falling/ in the context of a finite-rank $K =
\sum_{i=1}^n \lambda_i \phi_i \otimes \overline{\phi_i}$. Given disjoint $A_1,
\ldots, A_r \subseteq E$ and non-negative $k_1, \ldots, k_r$ summing to $k$,
it will be convenient to 
write $\kappa(j) := \min \big\{m \ge 1 \st j \le \sum_{\ell=1}^m k_\ell
\big\}$ for $j \le k$.
We have by Cauchy-Binet
\begin{align*}
\bE^K \Big[\prod_{\ell=1}^r &\big(\sX(A_\ell)\big)_{k_\ell}\Big]
=
\int_{\prod_{\ell=1}^r A_\ell^{k_\ell}} \rho_k \,d\mu^k
=
\int_{\prod_{\ell=1}^r A_\ell^{k_\ell}} \det (K \restrict (x_1, \ldots, x_k))
\,\prod_{j=1}^k d\mu(x_j)
\\
&=
\int_{\prod_{\ell=1}^r A_\ell^{k_\ell}}
\sum_{\substack{\sigma, \tau \in \Sym(k, n)
\\ \img (\sigma) = \img (\tau)}}
(-1)^\sigma (-1)^\tau \prod_{j=1}^k \lambda_{\sigma(j)} \phi_{\sigma(j)}(x_j)
\overline{\phi_{\tau(j)}(x_j)}
\,\prod_{j=1}^k d\mu(x_j)
\\
&=
\sum_{\substack{\sigma, \tau \in \Sym(k, n)
\\ \img (\sigma) = \img (\tau)}}
(-1)^\sigma (-1)^\tau 
\prod_{j=1}^k 
\int_{A_{\kappa(j)}}
\lambda_{\sigma(j)} \phi_{\sigma(j)}(x_j)
\overline{\phi_{\tau(j)}(x_j)}
\,d\mu(x_j)
\\
&=
\sum_{\substack{\sigma, \tau \in \Sym(k, n)
\\ \img (\sigma) = \img (\tau)}}
(-1)^\sigma (-1)^\tau 
\lambda^{\img(\sigma)}
\prod_{j=1}^k 
\bigip{\boI{A_{\kappa(j)}}
\phi_{\sigma(j)},
\overline{\phi_{\tau(j)}}}
\\
&=
\sum_{\sigma \in \Sym(k, n)}
(-1)^\sigma 
\lambda^{\img(\sigma)}
\det \Big[
\Bigip{\boI{A_{\kappa(j)}}
\phi_{\sigma(j)},
\overline{\phi_\ell}} \Big]_{\substack{j \le k \hfill \\ \ell \in \img(\sigma)}}
\,.
\end{align*}

As an immediate consequence of this formula, we obtain
the following important
principle of Goldman \cite[Proposition 12]{Goldman} 
that allows one to infer properties of continuous determinantal
point processes from corresponding properties of discrete determinantal
probability measures:

\procl t.transfer
Let $(E, \mu)$ and $(F, \nu)$ be two Radon measure spaces on locally
compact Polish sets.
Let $\Seq{A_i}$ be pairwise disjoint Borel subsets of $E$ and
$\Seq{B_i}$ be pairwise disjoint Borel subsets of $F$.
Let $\lambda_k \in [0, 1]$ with $\sum_k \lambda_k < \infty$.
Let $\Seq{\phi_k}$ be orthonormal in $L^2(E, \mu)$ and
$\Seq{\psi_k}$ be orthonormal in $\medstrut L^2(F, \nu)$.
Let $K := \sum_k \lambda_k \phi_k \otimes \overline{\phi_k}$
and $L := \sum_k \lambda_k \psi_k \otimes \overline{\psi_k}$.
If $\bigip{\boI{A_i} \phi_j, \phi_k} = \bigip{\boI{B_i} \psi_j, \psi_k}$
for all $i, j, k$, then the $\bP^K$-distribution of
$\Seq{\sX(A_i)}$ equals the $\bP^L$-distribution of
$\Seq{\sX(B_i)}$.
\endprocl

\rproof
When only finitely many $\lambda_k \ne 0$, this follows from our previous
calculation.
The general case follows from weak convergence of the processes
corresponding to the partial sums, as in the paragraph following \rref
l.weaklimit/.
%since the joint intensity functions
%converge and determine the count distributions via \rref e.disj/.
\Qed

This permits us to compare to discrete 
measures via \cite[Lemma 16]{Goldman}: 

\procl l.compare
Let $\mu$ be a Radon measure on a locally compact Polish space, $E$.
Let $\Seq{A_i}$ be pairwise disjoint Borel subsets of $E$.
Let $\phi_k \in L^2(E, \mu)$ for $k \ge 1$.
Then there exists a denumerable set $F$, pairwise disjoint subsets
$\Seq{B_i}$
of $F$, and $v_k \in \ell^2(F)$ such that
$\ip{\phi_j, \phi_k} = \ip{v_j, v_k}$ and
$\bigip{\boI{A_i} \phi_j, \phi_k} = \bigip{\boI{B_i} v_j, v_k}$
for all $i, j, k$.
\endprocl

\rproof
Without loss of generality, we may assume that $\bigcup_i A_i = E$.
For each $i$, fix an orthonormal
basis $\Seq{w_{i, j} \st j < n_i}$ for the subspace of $L^2(E, \mu)$ spanned by
$\{\boI{A_i} \phi_j\}$. Here, $n_i \in \bbN \cup \{\infty\}$.
Define $B_i := \{(i, j) \st j < n_i\}$ and $F := \bigcup_i B_i$.
Let $T$ be the isometric isomorphism from the span of $\{w_{i, j} \st i
\ge 1,\, j < n_i\}$ to $\ell^2(F)$ that sends $w_{i, j}$ to $\boI{\{(i,
j)\}}$.
Defining $v_k := T(\phi_k)$ yields the desired vectors.
\Qed

%[I think it follows from \cite[(7.6)]{Lyons:det} that there is always a trivial
%tail. One has to extend the above to locally trace-class $K$. Then look at
%events that depend only on finitely many counts. (7.6) gives control in
%terms of the number of counts needed to specify
%the first event, where the second is the tail event, to show
%asymptotic independence of the first event from the second, whence every
%tail event is trivial.]

\bsubsection{Stochastic Inequalities}{S.cinequalities}

We now show how the discrete models of \rref S.transf/ allow us to obtain
the analogues of the stochastic inequalities known to hold for discrete
determinantal probability measures. 
%Note that $E$ is second countable (i.e., has a countable basis),
%which implies that $E$ is $\sigma$-compact. 

For a Borel set $A
\subseteq E$, let $\sF(A)$ denote the $\sigma$-field 
on $\cN(E)$ generated by the functions $\xi \mapsto \xi(B)$ for Borel $B
\subseteq A$.
We say that a function that is measurable with respect to $\sF(A)$ is, more
simply, measurable with respect to $A$.
%A standard argument shows that if $f \colon \cN(E) \to \bbR$ is measurable
%with respect to $A$, then there are countably many $B_i \subseteq A$ such
%that $f$ is measurable with respect to the maps $\xi \mapsto \xi(B_i)$.
The obvious partial order on $\cN(E)$ allows us to define what it means for
a function $f \colon \cN(E) \to \bbR$ to be \dfnterm{increasing}.
As in the discrete case, we say that $\bP$ has \dfnterm{negative
associations} if
$
\bE[f_1 f_2] \le \bE[f_1] \bE[f_2] 
$
for every pair $f_1$, $f_2$ of bounded 
increasing functions that are measurable with
respect to complementary subsets of $E$.
An event is increasing if its indicator is increasing.
Then $\bP$ has negative associations iff
\rlabel{e.cnegass}
{
\bP(\ev A_1 \cap \ev A_2) \le \bP(\ev A_1) \bP(\ev A_2)
}
for every pair $\ev A_1$, $\ev A_2$ of 
increasing events that are measurable with
respect to complementary subsets of $E$.

We also say that \dfnterm{$\bP_1$ is stochastically dominated by $\bP_2$}
and write $\bP_1 \dom \bP_2$ if $\bP_1(\ev A) \le \bP_2(\ev A)$ for every
increasing event $\ev A$.

Call an event \dfnterm{elementary increasing} if it has the form $\{\xi \st
\xi(B) \ge k\}$, where $B$ is a relatively compact
Borel set and $k \in \bbN$.
Write $\isc(A)$ for the closure under
finite unions and intersections of the collection of
elementary increasing events with $B \subseteq A$; the notation
$\isc$ is chosen for ``upwardly closed".
Note that every event in $\isc(A)$ is measurable with respect to some
finite collection of functions
$\xi \mapsto \xi(B_i)$ for pairwise {\it disjoint\/} relatively compact
Borel $B_i \subseteq A$.
Write $\overline{\isc(A)}$ for the closure of $\isc(A)$ under monotone
limits, i.e., under unions of increasing sequences and under intersections
of decreasing sequences; these events are also increasing.
This is the same as the closure of $\isc(A)$ under countable unions and
intersections.

\procl l.approxincr
Let $A$ be a Borel subset of a locally compact Polish space, $E$.
Then $\overline{\isc(A)}$ is exactly the class of increasing Borel sets in
$\sF(A)$.
\endprocl

We give a proof at the end of this subsection.
First, we derive two consequences. 
A weaker version (negative correlations of elementary increasing events) of
the initial one is due to \rref b.Ghosh/.

\procl t.cFM
Let $\mu$ be a Radon measure on a locally compact Polish space, $E$.
Let $K$ be a locally trace-class positive contraction on $L^2(E, \mu)$.
Then $\bP^K$ has negative associations.
\endprocl

\rproof
Let $A \subset E$ be Borel.
Let $\ev A_1 \in \isc(A)$ and
$\ev A_2 \in \isc(E \setminus A)$.
Then $\ev A_1, \ev A_2 \in \sF(B)$ for some compact $B$ by definition
of $\isc(\cdot)$.
We claim that \rref e.cnegass/ holds for $\ev A_1$, $\ev A_2$,
and $\bP = \bP^{K_B}$, i.e., for $\bP = \bP^K$.
%\rlabel{e.1}
%{
%\bP^K(\ev A_1 \cap \ev A_2) \le \bP^K(\ev A_1) \bP^K(\ev A_2)
%\,.
%}

Now $\ev A_1$ is measurable with respect to a finite number of count
functions $\xi \mapsto \xi(B_i)$ for some disjoint $B_i \subseteq A \cap B$
($1 \le i \le n$) and
likewise $\ev A_2$ is measurable with respect to a finite number of 
functions $\xi \mapsto \xi(C_i)$ for some disjoint $C_i \subseteq B 
\setminus A$ ($1 \le i \le n$).
Thus, there are functions $g_1$ and $g_2$ such that $\boI{\ev A_1}(\xi) =
g_1\big(\xi(B_1), \dots, \xi(B_n)\big)$ and $\boI{\ev A_2}(\xi) =
g_2\big(\xi(C_1), \dots, \xi(C_n)\big)$.
By \rref t.transfer/ and
\rref l.compare/, there is some discrete determinantal probability
measure $\bP^Q$ on some denumerable set $F$ and pairwise disjoint sets
$B'_i, C'_i \subseteq F$ such that the joint $\bP^{K_B}$-distribution of all
$\sX(B_i)$ and $\sX(C_i)$ is equal to the joint $\bP^Q$-distribution of all
$\sX(B'_i)$ and $\sX(C'_i)$.
Define the corresponding events $\ev A'_i$ by
$\boI{\ev A'_1}(\xi) =
g_1\big(\xi(B'_1), \dots, \xi(B'_n)\big)$ and $\boI{\ev A'_2}(\xi) =
g_2\big(\xi(C'_1), \dots, \xi(C'_n)\big)$.
Since $\ev A'_i$
depend on disjoint subsets of $F$,
\rref t.FM/ gives that $\bP^Q(\ev A'_1 \cap \ev A'_2) \le \bP^Q(\ev
A'_1) \bP^Q(\ev A'_2)$. This is the same as \rref e.cnegass/
by \rref t.transfer/.
%and the fact that \rref e.1/ is really about $\bP^{K_B}$.

The same \rref e.cnegass/ clearly then holds in the less restrictive
setting $\ev A_i \in \overline{\isc(A)}$ by taking monotone limits.
\rref l.approxincr/ completes the proof.
\Qed

%The second consequence is Theorem 3 of \rref b.Goldman/:

\rprocl t.cdom {Theorem 3 of \rref b.Goldman/}
Suppose that $K_1$ and $K_2$ are two locally trace-class positive
contractions such that $K_1 \lLoew K_2$.
Then $\bP^{K_1} \dom \bP^{K_2}$.
\endprocl

\rproof
It suffices to show that $\bP^{K_1}(\ev A) \le \bP^{K_2}(\ev A)$ for every
$\ev A \in \isc(E)$.
Again, it suffices to assume that $K_i$ are trace class.
\rref l.compare/ applied
to all eigenfunctions of
$K_1$ and $K_2$ yields a denumerable $F$ and two positive
contractions $K'_i$ on $\ell^2(F)$, together with an event $\ev A'$, such
that $\bP^{K'_i}(\ev A') = \bP^{K_i}(\ev A)$ for $i = 1, 2$.
Furthermore, by construction,
every function in $\ell^2(F)$ is the image of a function in $L^2(E)$ under
the isometric isomorphism $T$ used to prove \rref l.compare/, whence
$K'_1 \lLoew K'_2$.
Therefore \rref t.dominate/ yields
$\smallstrut \bP^{K'_1}(\ev A') \le \bP^{K'_2}(\ev A')$, as desired.
\Qed

Again, it would be very interesting to have a natural monotone coupling of
$\bP^{K_1}$ with $\bP^{K_2}$. For some examples where this would be
desirable, see \rref S.orthogpoly/.

\rref l.approxincr/ will follow from this
folklore variant of a
theorem of Dyck \rref b.Dyck/:

\procl t.Dyck
Let $X$ be a Polish space on which $\le$ is a partial ordering that is
closed in $X \times X$. Let $\isc$ be a collection of
open increasing sets that generates the Borel subsets of $X$. Let $\isc^*$
be the closure of $\isc$ under countable intersections
and countable unions. Suppose that for all $x, y \in X$, either
$x \le y$ or there is $U \in \isc$ and an open set $V \subset X$
such that $x \in U$, $y \in V$, and $U \cap V = \emptyset$.
Then $\isc^*$ equals the class of increasing Borel sets.
\endprocl

\begin{proof}
Obviously every set in $\isc^*$ is Borel and increasing. To show the
converse, we prove a variant of Lusin's separation theorem.
Namely, we show that if $W_1 \subset X$ is increasing and
analytic (with respect to the
paving of closed sets, as usual) and if $W_2 \subset X$ is
analytic with $W_1 \cap W_2 = \emptyset$, then there exists $U \in \isc^*$ such
that $W_1 \subseteq U$ and $U \cap W_2 = \emptyset$. Taking $W_1$ to be
Borel and $W_2 := X \setminus W_1$ forces $U = W_1$ and gives the desired
conclusion.

To prove this separation property, we first show a stronger conclusion in
a special case: Suppose that $A_1, A_2
\subset X$ are compact such that $A_1$ is contained in an increasing
set $W_1$ that is disjoint from $A_2$; then there exists an open
$U \in \isc^*$ and an open $V$ such
that $A_1 \subseteq U$, $A_2 \subseteq V$, and $U \cap V = \emptyset$. 
Indeed, since $W_1$ is increasing, for every $(x, y) \in A_1 \times A_2$,
we do {\it not} have that $x \le y$, whence by hypothesis, there exist
$U_{x, y} \in \isc$ and an open $V_{x, y}$ with $x \in U_{x, y}$, $y \in
V_{x, y}$, and $U_{x, y} \cap V_{x, y} = \emptyset$.
Because $A_2$ is compact, for each $x \in A_1$, we may choose
$y_1, \ldots, y_n \in A_2$ such that $A_2 \subseteq V_x :=
\bigcup_{i=1}^n V_{x, y_i}$. Define 
$U_x := \bigcap_{i=1}^n U_{x, y_i}$.
Then $U_x$ is open, contains $x$, and is disjoint from $V_x$, whence
compactness of $A_1$ ensures the existence of $x_1, \ldots, x_m \in
A_1$ with $A_1 \subseteq U := \bigcup_{j=1}^m U_{x_j} \in
\isc^*$. Then $V := \bigcap_{j=1}^m V_{x_j}$ is open, contains $A_2$,
and is disjoint from $U$, as desired.

To prove the general case, let $\pi_1$ and $\pi_2$ be the two coordinate
projections on $X^2 = X \times X$. Define $I(A) =
I\big(\pi_1(A) \times \pi_2(A)\big)$ for $A \subseteq X^2$ to
be 0 if there exists $U \in \isc^*$ such that $\pi_1(A) \subseteq U$ and $U
\cap \pi_2(A) = \emptyset$; and to be 1 otherwise. 
%Clearly $I(A) = I\big(\pi_1(A) \times \pi_2(A)\big)$.

We claim that $I$ is a capacity in the sense of \cite[(30.1)]{Kechris}.
It is obvious that $I(A) \le I(B)$ if $A \subseteq B$ and it is simple to
check that if $A_1 \subseteq A_2 \subseteq \cdots$, then $\lim_{n
\to\infty} I(A_n) = I\big(\bigcup_n A_n\big)$.
Suppose for the final property that $A$ is compact and $I(A) = 0$; we must
find an open $B \supseteq A$ for which $I(B) = 0$.
%If $A = \emptyset$, then we may take $B = A$.
There exists some $W_1 \in \isc^*$ with $\pi_1(A) \subseteq W_1$ and
$W_1 \cap \pi_2(A) = \emptyset$.
Then the result of the second paragraph yields sets $U$ and $V$ that
give $B := U \times V$ as desired.

Now let $W_1$ and $W_2$ be as in the first paragraph.
If $A \subseteq W_1 \times W_2$ is compact,
then setting $A_i := \pi_i(A)$ and applying the second paragraph shows that
$I(A) = 0$.
Thus, by the Choquet capacitability theorem \cite[(30.13)]{Kechris}, $I(W_1
\times W_2) = 0$.
\end{proof}

\begin{proofof}{\rref l.approxincr/}
Clearly every set in $\overline{\isc(A)}$ is increasing and in $\sF(A)$.
For the converse, endow $A$ with a metric so that it becomes locally compact
Polish while preserving its class of relatively compact sets and its Borel
$\sigma$-field: Choose a denumerable partition of $A$ into relatively compact
sets $A_i$ and make each one compact and of diameter at most 1; make the
distance between $x$ and $y$ be 1 if $x$ and $y$ belong to different $A_i$.
Let $X := \cN(A)$ with the vague topology and let $\isc$ be the class of
elementary increasing events defined with respect to (relatively compact)
sets $B \subseteq A$ that are open for this new metric.
Apply \rref t.Dyck/. Since
$\isc^* \subseteq \overline{\isc(A)}$, the result follows.
\end{proofof}

\bsubsection{Example: Orthogonal Polynomial Ensembles}{S.orthogpoly}

Natural examples of determinantal point processes arise from
orthogonal polynomials with respect to a probability measure $\mu$ on $\bbC$.
Assume that $\mu$ has infinite support and finite moments of all orders.
Let $K_n$ denote the orthogonal projection of $L^2(\bbC, \mu)$ onto the
linear span $\poly_n$ of the functions $\{1, z, z^2, \ldots, z^{n-1}\}$.
%i.e., onto the polynomials of degree less than $n$.
There exist unique (up to signum) polynomials $\phi_k$ of degree $k$ such
that for every $n$,
$\Seq{\phi_k \st 0 \le k < n}$ is an orthonormal basis of $\poly_n$.
By elementary row operations, we see that for variables $(z_1, \ldots,
z_n)$, the map $(z_1, \ldots, z_n) \mapsto \det [\phi_i(z_j)]_{i, j \le n}$
is a Vandermonde
polynomial up to a constant factor, whence 
$$
\det (K_n \restrict \{z_1, \ldots, z_n\}) 
=
\det [\phi_i(z_j)] [\phi_i(z_j)]^*
=
c_n \prod_{1 \le i < j \le n} |z_i-z_j|^2
$$
for some constant $c_n$.
Therefore, the density of $\bP^{K_n}$ (with points randomly ordered)
with respect to $\mu^n$ is given by
$c_n/n!$ times the square of a Vandermonde determinant.

Classical examples include the following:
\begin{enumerate}[label=\bfseries OPE\arabic*.]

\item
If $\mu$ is Gaussian measure on $\bbR$, i.e., $d\mu(x) =
(2\pi)^{-1/2} e^{-x^2/2} \,dx$, then $\phi_k$ are the Hermite polynomials,
$c_n = \big(\prod_{j=1}^{n-1} j!\big)^{-1}$, and $\bP^{K_n}$ is the law of the
\dfnterm{Gaussian unitary ensemble}, which is the set of eigenvalues of
$(\rM+\rM^*)/\sqrt2$, where $\rM$ is an $n \times n$ matrix whose entries are
independent standard complex Gaussian. (A standard complex Gaussian random
variable is the same as a standard Gaussian vector in $\bbR^2$ divided by
$\sqrt 2$ in order that the complex variance equal 1. Its density
is $\pi^{-1} e^{-|z|^2}$ with respect to Lebesgue measure on $\bbC$.) This is
due to Wigner; see \rref b.Mehta/.

\item 
If $\mu$ is unit Lebesgue measure on the unit circle $\{z \st |z| = 1\}$,
then $\phi_k(z) = z^k$, so $c_n = 1$, and $\bP^{K_n}$ is the law of the
\dfnterm{circular unitary ensemble}, which is the set of eigenvalues of a
random matrix whose distribution is Haar measure on the set of $n \times n$
unitary matrices. This ensemble was introduced by Dyson, but the law of the
eigenvalues is due to Weyl; see \rref b.HKPV:book/.

\item
If $\mu$ is standard Gaussian measure on $\bbC$, then $\phi_k(z) =
z^k/\sqrt{k!}$, $c_n = \big(\prod_{j=1}^{n-1} j!\big)^{-1}$, and $\bP^{K_n}$ is
the law of the \dfnterm{$n$th (complex) Ginibre process}, which is the
set of eigenvalues of an $n \times n$ matrix whose entries are independent
standard complex Gaussian. This is due to Ginibre; see \rref
b.HKPV:book/.

\item
If $\mu$ is unit Lebesgue measure on the unit disk $\disk := \{z \st |z| < 1\}$,
then $\phi_k(z) = \sqrt{k+1}\, z^k$, so $c_n = n!$, and the limit of
$\bP^{K_n}$ is the law of the zero set of the random power series whose
coefficients are independent standard complex Gaussian, which converges in
the unit disk a.s. This is due to Peres and Vir\'ag \rref b.PeresVirag/.

\item
If $\mu$ has density $z \mapsto n \pi^{-1}\big(1+|z|^2\big)^{-n-1}$ with
respect to Lebesgue measure on $\bbC$, then $\phi_k(z) = \sqrt{\binom{n-1}
{k} } z^k$ for $k < n$, so $c_n = \prod_{j=1}^{n-1} \binom{n-1} {j}$,
and $\bP^{K_n}$ is
the law of the \dfnterm{$n$th spherical ensemble}, which is the set of
eigenvalues of $\rM_1^{-1} \rM_2$ when $\rM_i$ are
independent $n \times n$ matrices whose entries are independent
standard complex Gaussian. (Here, we are limited to $\poly_n$ since the
larger spaces do not lie in $L^2(\mu)$.) This is due to Krishnapur
\rref b.Krishnapur:thesis/; see
\rref b.HKPV:book/. The process was studied earlier by \rref b.Caillol/ and
\rref b.FJM/, but without observing the connection to eigenvalues.
Inverting stereographic projection, we identify this process with one
whose density with respect to Lebesgue measure on 
the unit sphere in $\bbR^3$
is proportional to
$\prod_{1 \le i < j \le n} \|\bv_i - \bv_j\|^2$.

\end{enumerate}

For additional information on such processes, see
\cite{Simon:CD,Hardy:average,RiderVirag:H1,Forrester:book}. 
For an extension to complex manifolds, see
\cite{Berman:bulk,Berman:DPP,Berman:sharp}.

By \rref t.cdom/, the processes $\bP^{K_n}$ stochastically increase in $n$
for each of the examples above except the last. It would be interesting to
see natural monotone couplings. Perhaps the
last example also increases stochastically in $n$.

The \dfnterm{Ginibre process} is the limit of the $n$th Ginibre
processes as $n \to\infty$; it has the kernel
$e^{z \bar w}$ with respect to standard Gaussian measure 
on $\bbC$.  This process is invariant under all isometries of $\bbC$.
% Its kernel with respect to Lebesgue measure is $\pi^{-1} e^{-|z-w|^2/2}$
% up to a factor $e^{i \Im z \bar w}$, which changes by translation via
% a gauge transformation and which is clearly invariant under rotation; it
% is also the same measure under reflection in the real axis
For each of the plane, sphere, and hyperbolic disk, there is only a 
1-parameter family of determinantal point processes having a kernel $K(z, w)$
that is holomorphic in $z$ and in $\bar w$ and whose law is 
isometry invariant \cite[Theorem 3.0.5]{Krishnapur:thesis}.
For the sphere, that family has already been given above; the parameter is
a positive integer. For the other two families, the parameter is a positive
real number, $\alpha$. In the case of the plane, the processes are 
related simply by homotheties, $M_\alpha \colon z \mapsto z/\alpha$. The
push-forward of the Ginibre process 
with respect to $M_{\sqrt\alpha}$ has kernel
$e^{\alpha z \bar w}$ with respect to the measure $\alpha \pi^{-1}
e^{-\alpha |z|^2} d\mu(z)$, where $\mu$ is Lebesgue measure on $\bbC$. Do
these processes increase stochastically in $\alpha$, like Poisson processes
do?
In the hyperbolic disk, the
processes have kernel $\alpha (1 - z \bar w)^{-\alpha-1}$
with respect to the measure $\pi^{-1} (1 - |z|^2)^{\alpha-1} d\mu(z)$,
where $\mu$ is Lebesgue measure on $\disk$. (We fix a branch of
$(1-z)^{-\alpha-1}$ for $z \in \disk$.) These give orthogonal
projections onto the generalized Bergman spaces. The case $\alpha = 1$ is that
of the limiting OPE4 above. Do these processes stochastically
increase in $\alpha$?

\bsection{Completeness}{s.complete}

Recall that when $H$ is a finite-dimensional subspace of $\eH$, the measure
$\bP^H$ is supported by those
subsets $B \subseteq E$ that project to a basis of $H$ under $P_H$.
Similarly, when $K$ is the kernel of a finite-rank orthogonal projection
onto $H \subset L^2(E, \mu)$, 
define the functions $K_x := K(\cdot, x) = \sum_{k \ge 1}
\overline{\phi_k(x)} \phi_k \in H$.
Then the measure
$\bP^K$ is supported by those $\xi$ such that $\Seq{K_x \st x \in \xi}$ is
a basis of $H$, since $K(x, y) = \ip{K_y, K_x}$.
Here, $x \in \xi$ means that $\xi(\{x\}) = 1$.

The question of extending this to infinite-dimensional $H$ turns out to be
very interesting. 
A basis of a finite-dimensional vector space is a minimal spanning set. 
Although $P_H \ba$ is $\bP^H$-a.s.\ linearly independent,
minimality does not hold in general, even for
the wired spanning forest of a tree, as shown by the examples in \rref
b.HeicklenLyons/.
See also \rref c.ell2min/.
% Perhaps invariance would make it hold. Then we would get bilaterally
% deterministic Bernoulli shifts. But I think not.
However, the other half of being a basis does hold in the discrete case and
is open in the continuous case.

\bsubsection{Discrete Completeness}{S.dcompl}

Let $[V]$ be the closed linear span of $V \subseteq \eH$.

\rprocl t.basis {\rref b.Lyons:det/}
For every $H \subsp \eH$, we have $[P_H \ba] = H$ $\bP^H$-a.s.
\endprocl

We give an application of \rref t.basis/ for $E = \bbZ$, 
but it has an analogous statement for every countable abelian group.
Let $\bbT := \bbR/\bbZ$ be the unit circle equipped with unit Lebesgue
measure. For a measurable function $f \colon \bbT \to \bbC$ and $n \in
\bbZ$, the \dfnterm{Fourier coefficient} of $f$ at $n$ is 
$\medstrut
\widehat f(n) := \int_\bbT f(t) e^{-2\pi i n t} \,dt
$.
Let $\widehat f \restrict S$ denote the restriction of $\widehat f$ to
$S$.
If $A \subseteq \bbT$ is measurable, we say $S \subseteq \bbZ$ is
\dfnterm{complete for} $A$ if the set $ \{ f \boI A \st f \in L^2(\bbT),\,
\medstrut \widehat f
\restrict (\bbZ \setminus S) \equiv 0 \}$ is dense in $L^2(A)$, 
where we identify $L^2(A)$ with
the set of functions in $L^2(\bbT)$ that vanish outside
$A$. 
%(This is the
%concept dual to restriction set, with the roles of $A$ and $S$ reversed. It is
%also the same as a uniqueness set for $A$, meaning that if $f \in L^2(A)$
%satisfies $\widehat f \restrict S \equiv 0$, then $f \equiv 0$.)
The case where $A$ is an interval is quite classical; see
\rref b.Redheffer/ for a review.
A crucial role in that case is played by the following notion of density of
$S$.

\procl d.BM
For an interval $[a, b] \subset \bbR \setminus \{ 0 \}$, define its
\dfnterm{aspect} 
$$
\asp\big([a, b]\big)
:=
\max \big\{ |a|, |b| \big\}/ \min \big\{ |a|, |b| \big\}
\,.
$$
For a discrete
$S \subseteq \bbR$, the \dfnterm{Beurling-Malliavin density} of $S$, denoted
$\bm(S)$, is the supremum of those $D \ge 0$ for which there exist disjoint
nonempty intervals $I_n \subset \bbR \setminus \{ 0 \}$ 
with $|S \cap I_n| \ge D |I_n|$ for all $n$ and $\sum_{n \ge 1}
[\asp(I_n)-1]^2 = \infty$.
\endprocl

\iffalse
A simpler form of the Beurling-Malliavin density was provided by \rref
b.Red:Two/, who showed that 
\rlabel e.BMRed
{
\bm(S) = \inf \left\{ c \st \exists \hbox{ an injection }\beta : S \to \bbZ
\hbox{ with } \sum_{k \in S} \Big|\frac{1}{k} - \frac{c}{\beta(k)}\Big| <\infty
\right\} 
\,.
}
\fi

\rprocl c.seqdual {\rref b.Lyons:det/}
Let $A \subset \bbT$ be Lebesgue measurable with measure $|A|$.
Then there is a set of Beurling-Malliavin density $|A|$ in $\bbZ$ that is 
complete for $A$.
Indeed,
let\/ $\bP^A$ be the determinantal probability measure on $2^\bbZ$
corresponding to the Toeplitz matrix $(j, k) \mapsto \widehat{\boI{A}}(k-j)$.
Then $\bP^A$-a.e.\ $S \subset \bbZ$ is complete for $A$ and has
$\bm(S) = |A|$.
\endprocl

When $A$ is an interval, the celebrated theorem of Beurling and Malliavin
\rref b.BM/ says that if
$S$ is complete for $A$, then $\bm(S) \ge |A|$, and that if
$\bm(S) > |A|$, then $S$ is complete for $A$.
(This holds for $S$ that
are not necessarily sets of integers, but we are concerned in this subsection
only with $S \subseteq \bbZ$.) 

\rref c.seqdual/ can be compared (take $\bbT \setminus A$ and $\bbZ
\setminus S$) to a theorem of Bourgain and Tzafriri
\rref b.BTz/, according to
which there is a set $S \subset \bbZ$ of (Schnirelman) density at least
$2^{-8} |A|$ such that if $f \in L^2(\bbT)$ and $\widehat f$ vanishes off
$S$, then 
$$
|A|^{-1} \int_A |f(t)|^2 \,dt \ge 2^{-16} \|f\|_2^2
\,.
$$
It would be interesting to find a quantitative strengthening of \rref
c.seqdual/ that would encompass
this theorem of \cite{BTz}.

The following theorem is equivalent to \rref t.basis/ by duality:

\rprocl t.morris {\rref b.Lyons:det/}
For every $H \subsp \eH$, we have $\overline{P_{[\ba]} H} =
[\ba]$ $\bP^H$-a.s.
\endprocl

As an example, consider the wired spanning forest of a graph,
$G$. Here, $H := \STAR(G)$.
In this case, $H_B := \overline{P_{[B]} \STAR(G)} = \STAR(B)$ for $B
\subseteq E$.
Thus, the conclusion of \rref t.morris/ is that $\bP^{H_\fo}$, which equals
$\wsf_\fo$, is concentrated on the singleton $\{\fo\}$ for
$\wsf_G$-a.e.\ $\fo$. 
This was a 
conjecture of \BLPSusf, established by %Morris
\rref b.Morris/.

\procl c.ell2min
For every $H \subsp \eH$, 
$\bP^H$-a.s.\ the maps $P_H \colon [\ba] \to H$ and $P_{[\ba]} \colon H \to
[\ba]$ are injective with dense image.
\endprocl

\rproof
Both statements are equivalent to $[\ba] \cap H^\perp = \{0\} = H \cap
\ba^\perp$, and these are the contents of Theorems \ref{thm:basis} and
\ref{thm:morris}.
\Qed

\iffalse
proved that on any network $(G, w)$ (where $G$ is the underlying
graph and $w$ is the function assigning conductances, or weights, to the
edges), for $\wsf(G, w)$-a.e.\ forest $\fo$ and for every component tree $T$
of $\fo$, the $\wsf$ of $(T, w \restrict T)$ equals $T$ a.s.
This suggested \rref b.Lyons:det/ the following extension. Given a
subspace $H$ of $\eH$ and a set $B \subseteq E$, the subspace of $[B]$
``most like" or ``closest to" $H$ is the closure of the image of $H$ under the
orthogonal projection $P_{[B]}$; we denote this
subspace by $H_B$. For example, if $H = \STAR(G)$,
then $H_B = \STAR(B)$ since
for each $x \in V(G)$, we have $P_{[B]} (\star_x^G) = \star_x^B$.
To say that $\bP^{H_B}$ is concentrated on $ \{ B \}$ is the same as to say
that $H_B = [B]$. This motivated the following theorem 
and shows how it is
an extension of Morris's theorem.
\fi

\bsubsection{Continuous Completeness}{S.ccompl}

If $K$ is a locally trace-class orthogonal projection onto $H$, then for $h
\in  H$, we have 
$$
h(x)
=
(Kh)(x)
=
\int_E K(x, y) h(y) \,d\mu(y)
=
\int_E h(y) \overline{K(y, x)} \,d\mu(y)
=
\bigip{h, K_x}
\,.
$$
In other words, $K$ is a reproducing kernel for $H$.
A subset $S$ of $H$ is called \dfnterm{complete for $H$} if the closed
linear span of $S$ equals $H$; equivalently, the only element of $H$ that
is orthogonal to $S$ is 0.

An analogue of \rref t.basis/ was conjectured by Lyons and Peres in 2010:

\procl g.cbasis
If $K$ is a locally trace-class orthogonal projection onto $H$, then
for $\bP^K$-a.e.\ $\sX$,
$[\{K_x \st x \in \sX\}] = H$, i.e.,
if $h \in H$ and $h \restrict \sX = 0$, then $h \equiv 0$.
\endprocl

Just as in the discrete case, this appears to be on the critical border for
many special instances, as we illustrate for several processes where $E =
\bbC$:
\begin{enumerate}[label=\bfseries \arabic*.]

\item
Let $\mu$ be Lebesgue measure on $\bbR$ and $K(x, y) := \sin
\pi(x-y)/\big(\pi(x-y)\big)$, the \dfnterm{sine-kernel process}. Denote the
Fourier transform on $\bbR$ by 
$
\widehat f(t)
:=
\int_\bbR f(x) e^{-2\pi i t x} \,dx
$
for $f \in L^1(\bbR)$, and, by isometric extension, for $f \in L^2(\bbR)$.
Write $I := \boI{[-1/2, 1/2]}$.
Since $K(x, 0) = \widehat{I}(x)$, 
we have 
$
(Kf)(x)
=
\big(f * \widehat I\big)(x)
=
\widehat{\check f I}(x)
$,
where $\check f$ is the inverse Fourier transform of $f$.
Therefore, the induced operator $K$
arises from the orthogonal projection onto the Paley-Wiener space
$\{f \in L^2(\bbR, \mu) \st \check f(t) = 0 \mbox{ if } |t| >
1/2\}$. The sine-kernel
process arises frequently; e.g., it is various scaling limits of
the $n$th Gaussian unitary ensemble ``in the bulk" as $n \to\infty$.
(A related scaling limit of the GUE is Wigner's semicircle distribution.)
We may more easily interpret \rref g.cbasis/ for Fourier transforms of
functions in $L^2[-1/2, 1/2]$: It says that for $\bP^K$-a.e.\ $\sX$, the
only $h \in L^2[-1/2, 1/2]$ such that $\widehat h \restrict \sX = 0$ is $h
\equiv 0$.
Although the Beurling-Malliavin theorem applies, 
no information can be deduced
because $\bm(\sX) = 1$ a.s.
However, Ghosh \rref b.Ghosh/ has proved this case.

\item
Let $\mu$ be standard Gaussian measure on $\bbC$ and $K(z, w) := e^{z \bar
w}$. This is the Ginibre process.
It corresponds to orthogonal projection onto the \dfnterm{Bargmann-Fock
space} $B^2(\bbC)$ consisting of the entire functions that lie in
$L^2(\bbC, \mu)$; this is the space of power series $\sum_{n \ge 0} a_n
z^n$ such that $\sum_n n! |a_n|^2 < \infty$.
Completeness of a set of elements $\medstrut\big\{e^{\lambda z} \st \lambda \in
\Lambda\big\} \subset B^2(\bbC)$ in $B^2(\bbC)$ is equivalent 
to completeness in $L^2(\bbR)$ (with Lebesgue measure) of the Gabor system
of windowed complex exponentials
$$
\Big\{t \mapsto \exp \big[-i \myIm \lambda t - (t -
\myRe \lambda)^2\big] \st \lambda \in \sqrt{2}\Lambda\Big\}
\,,
$$
which is used in time-frequency
analysis of non-band-limited signals. The equivalence is proved using the
Bargmann transform 
$$
f \mapsto
\Big(z \mapsto
\pi^{-1/4} \int_\bbR f(t) \exp \big[\sqrt{2} t z -
\frac{z^2}{2} - \frac{t^2}{2}\big] \,dt \Big)
\,,
$$
which is an isometry from $L^2(\bbR)$ to $B^2(\bbC)$.
That the critical density is 1 was shown in various senses going back to
von Neumann; see \rref b.CLP/.
This case has also been proved by Ghosh \rref b.Ghosh/.

\item
Let $\mu$ be unit 
Lebesgue measure on the unit disk $\disk := \{z \st |z| <
1\}$ and $K(z, w) := \big(1 - z \bar w\big)^{-2}$.
This process is the limiting OPE4 in \rref S.orthogpoly/.
%discussed earlier, namely, it is
%the law of the zero set of the random power series whose
%coefficients are independent standard complex Gaussian.
It corresponds to orthogonal projection onto the \dfnterm{Bergman space}
$A^2(\disk)$ consisting of the analytic functions that lie in $L^2(\disk,
\mu)$.
What is known about the zero sets of functions in the Bergman space
\rref b.Duren/ is insufficient to
settle \rref g.cbasis/ in this case and it remains open.

\end{enumerate}

The two instances above that have been proved by Ghosh \rref b.Ghosh/ follow
from his more general result that \rref g.cbasis/ holds whenever $\mu$ is
continuous and $\bP^K$ is \dfnterm{rigid}, which means that $\sX(B)$ is
measurable with respect to the $\bP^K$-completion of $\sF(E \setminus B)$
for every ball $B \subset E$.
The limiting process OPE4 is not rigid \rref b.HS:tolerance/.
Ghosh and Krishnapur (personal communication, 2014) have shown that $\bP^K$
is rigid only if $K$ is an orthogonal projection.
It is not sufficient that $K$ be a projection, as the example of the
Bergman space shows.
A necessary and sufficient condition to be rigid is not known.

Let $K$ be a locally trace-class orthogonal projection onto $H \subsp
L^2(E, \mu)$.
For a function $f$, write $f_K$ for the function $f(x)/\sqrt{K(x, x)}$.
Let $\sX \sim \bP^K$. Clearly 
$f_K \restrict \sX \in \ell^2(\sX)$ for a.e.\ $\sX$. Also, for $h \in
H$, the function $h_K$ is bounded. A conjecture analogous to \rref
c.ell2min/ is that $\sX$ is a sort of set of interpolation for $H$ in the
sense that given any countable dense set $H_0 \subset H$, for a.e.\ $\sX$, the
set $\{h_K \restrict \sX \st h \in H_0\}$ is dense in $\ell^2(\sX)$.

One may also ask about completeness for appropriate Poisson point processes.

\bsection{Discrete Invariance}{s.dinvar}

Suppose $\gp$ is a group that acts on $E$ and that $K$ is $\gp$-invariant,
i.e., $K(\gpe x, \gpe y) = K(x, y)$ for all $\gpe \in \gp$, $x \in E$, and
$y \in E$. (This is equivalent to the operator $K$ being
$\gp$-equivariant.) Then the probability measure $\bP^K$ is $\gp$-invariant.
This contact with ergodic theory and other areas of mathematics
suggests many interesting questions.
Lack of
space prevents us from considering more than just a few aspects of the case
where $E$ is discrete and from giving all definitions.

\bsubsection{Integer Lattices}{S.zd}

Let $E := \gp := \bbZ^d$.
In this case, $K$ is invariant iff $\medstrut 
K(m, n) = \widehat f(n-m)$ for some $f
\colon \bbT^d \to [0, 1]$, where 
$\widehat f(n) := \int_{\bbT^d} f(t) e^{-2\pi i n \cdot t} \,dt$.
We write $\bP^f$ in place of $\bP^K$.
Some results and questions from \rref b.LS:dyn/ follow.

\procl t.Bern
For all $f$, the process $\bP^f$ is isomorphic to a Bernoulli process. 
\endprocl

This was shown in dimension 1
by \rref b.ShiTak:II/ for those $f$ such that $\sum_{n \ge 1}
n |\widehat f(n)|^2 <\infty$ by showing that those $\bP^f$ are weak
Bernoulli (WB),
also called ``$\beta$-mixing" and ``absolutely
regular". Despite its name, it is known that WB is
strictly stronger than Bernoullicity. 
The precise class of $f$ for which $\bP^f$ is WB is not known.

As usual, the \dfnterm{geometric mean} of a nonnegative function $f$ is
$
\GM(f) := \exp \int \log f
$.

\procl t.GMdom
For all $f$, the process $\bP^f$ stochastically dominates product measure
$\bP^{\GM(f)}$ and is stochastically dominated by product measure
$\bP^{1-\GM(\bfone-f)}$. These bounds are optimal.
\endprocl

We conjecture that (Kolmogorov-Sinai) entropy is concave,
as would follow from \rref g.concave/.

\procl g.invconcave
For all $f$ and $g$, we have $H\big(\bP^{(f+g)/2}\big) \ge \big(H(\bP^f) +
H(\bP^g)\big)/2$.
\endprocl

\procl q.block
Let $f \colon \bbT \to [0, 1]$ be a trigonometric polynomial of
degree $m$.
Then $\bP^f$ is $m$-dependent, as are all $(m+1)$-block factors of
independent processes.
Is $\bP^{f}$ an $(m+1)$-block factor of an i.i.d.\ process?
This is known when $m = 1$ \rref b.Broman/.
\endprocl

\bsubsection{Sofic Groups}{S.sofic}

Let $\gp$ be a sofic group, a class of groups that includes all finitely
generated amenable
groups and all finitely generated residually amenable groups. No finitely
generated group is known not to be sofic.
Let $E$ be $\gp$ or, more generally, a set acted on by $\gp$ with finitely
many orbits, such as the edges of a Cayley graph of $\gp$.
The following theorems are from \rref b.LyonsThom/.

\procl t.sBern
For every $\gp$-equivariant positive contraction $Q$ on $\eH$, the process
$\bP^Q$ is a $\dbar$-limit of finitely dependent (invariant) processes.
If\/ $\gp$ is amenable and $E = \gp$, then $\bP^Q$ is isomorphic to a Bernoulli
process. 
\endprocl

Even if $\bP^1$ and $\bP^2$ are $\gp$-invariant probability measures
on $2^\gp$ with $\bP^1 \dom \bP^2$, there need not be a
$\gp$-invariant monotone coupling of $\bP^1$ and $\bP^2$ \rref
b.Mester:mono/.
The proof of the preceding theorem depends on the next one:

\procl t.monojoin
If\/ $Q_1$ and $Q_2$ are two $\gp$-equivariant positive
contractions on $\eH$ with $Q_1 \lLoew Q_2$,
then there exists a $\gp$-invariant monotone coupling 
of\/ $\bP^{Q_1}$ and\/ $\bP^{Q_2}$.
\endprocl

The proof of \rref t.sBern/ also uses the inequality
$$
\dbar\big(\bP^Q, \bP^{Q'}\big) \le 6 \cdot 3^{2/3} \|Q - Q'\|_1^{1/3}
$$
for equivariant positive contractions, $Q$ and $Q'$, where $\|T\|_1 := \tr
(T^* T)^{1/2}$ is
the Schatten 1-norm.
When $Q$ and $Q'$ commute,
one can improve this bound to
$$\dbar(\bP^Q, \bP^{Q'}) \le \|Q - Q'\|_1\,.$$
We do not know
whether this inequality always holds.

Write $\FK(Q) := \exp \tr \log |Q|$ for the Fuglede-Kadison determinant of
$Q$ when $Q$ is a $\gp$-equi\-var\-iant operator. The following would
extend \rref t.GMdom/. It is open even for finite groups.

\procl g.FKdom
For all\/ $\gp$-equivariant positive contractions $Q$ on $\ell^2(\gp)$,
the process $\bP^Q$ stochastically dominates product measure
$\bP^{\FK(Q)I}$ and is stochastically dominated by product measure
$\bP^{I-\FK(I-Q)I}$, and these bounds are optimal.
\endprocl

\bsubsection{Isoperimetry, Cost, and $\ell^2$-Betti Numbers}{S.betti}

It turns out that the expected degree of a vertex in the free uniform
spanning forest of a Cayley graph depends only on the group, via its first
$\ell^2$-Betti number, $\beta_1(\gp)$, and not on the
generating set used to define the Cayley graph \rref b.Lyons:betti/:

\procl t.betti In every Cayley graph $G$ of a group $\gp$, we have 
$$
\bE_{\fsf(G)}[\deg_\fo (\bp)] = 2 \beta_1(\gp) + 2
\,.
$$
\endprocl

This is proved using the representation of $\fsf$ as a determinantal
probability measure.
It can be used to give a uniform bound on expansion constants \rref
b.LPV/:

\procl t.LPV For every finite symmetric
generating set $S$ of a group $\gp$, we have 
$
{|S A \setminus A|}
>
2 \beta_1(\gp) {|A|}
$
for all finite non-empty $A \subset \gp$.
\endprocl

There are extensions of these results to higher-dimensional CW-complexes
and higher $\ell^2$-Betti numbers \rref b.Lyons:betti/.

In unpublished work with D.~Gaboriau \cite{LG:approach}, we have shown the following:

\procl t.damien
Let $G$ be a Cayley graph of a finitely generated group $\gp$ and $\epsilon
> 0$. Then there exists a $\gp$-invariant finitely dependent determinantal
probability measure $\bP^Q$ on $\{0, 1\}^{\edges(G)}$ that stochastically
dominates $\fsf_G$ and such that 
$$
\bE^Q\big[\deg_\qba(\bp)\big] \le
\bE_\fsf\big[\deg_\fo(\bp)\big] + \epsilon
\,.
$$
In addition, if\/ $\gp$ is sofic, then $\dbar(\bP^Q, \fsf) \le \epsilon$.
\endprocl

If it could be shown that $\bP^Q$, or indeed every invariant
finitely dependent probability measure that dominates $\fsf$, yields a
connected subgraph
a.s., then it would follow that $\beta_1(\gp) + 1$ is equal to the cost of
$\gp$, a major open problem of \rref b.Gaboriau:invar/.

%\bsubsection{Higher-Dimensional Complexes}{S.top}

%\bsection{Continuous Invariance}{s.cinvar}
%
%\bsubsection{Complex Plane}{S.complex}
%
%\bsubsection{Hyperbolic Plane}{S.hyperb}
%
%\bsubsection{Sphere}{S.sphere}

\section*{Acknowledgments}
I am grateful to Alekos Kechris for informing me of \rref t.Dyck/;
the proof given seems to be due to Alain Louveau.
I thank Norm Levenberg for references.

%\cleardoublepage

\phantomsection

\addcontentsline{toc}{section}{References}

%\frenchspacing

\end{document}